\documentclass[10pt, a4paper]{article}
\usepackage{latexsym}
\usepackage{amsmath}
\usepackage{amssymb}
\usepackage[all]{xy}
\usepackage{bm}
\usepackage{amsfonts}
\usepackage{verbatim}
\usepackage{amsthm}
\usepackage{mathrsfs}
\usepackage{epsfig}
\usepackage{xy}
\usepackage{array}
\usepackage{stmaryrd}
\usepackage{graphicx,color}
\usepackage{xcolor}
\usepackage{tikz}
\usetikzlibrary{arrows,calc}
\usepackage{etex}
\usepackage{mathdots}
\usepackage{float}
\usepackage{graphics}
\usepackage{pdflscape}
\usepackage{cite}

\usepackage{anysize,hyperref}
\input xypic
\xyoption{all}

\usepackage[perpage,symbol]{footmisc}
\usepackage{setspace}
\renewcommand{\cal}{\mathcal}
\def\A{\mathscr{A}}
\def\B{\mathscr{B}}
\def\C{\mathscr{C}}
\def\E{\mathbb{E}}
\def\s{\mathfrak{s}}

\def\op{^\mathrm{op}}
\def\Ab{\mathit{Ab}}
\def\del{\delta}
\def\dr{\ar@{->}[r]}
\def\X{\mathscr{X}}
\def\Y{\mathscr{Y}}
\def\Z{\mathscr{Z}}

\def\X{\mathscr{X}}
\def\Y{\mathscr{Y}}

\def\Hom{\mbox{Hom}}

\begin{document}

\baselineskip=15pt
\title{\Large{\bf Tilting subcategories in extriangulated categories}\footnotetext{This work was supported by the NSF of China (Grants No.11671221)}}
\medskip
\author{\textbf{Bin Zhu and Xiao Zhuang}}

\date{}

\maketitle
\def\blue{\color{blue}}
\def\red{\color{red}}

\newtheorem{theorem}{Theorem}[section]
\newtheorem{lemma}[theorem]{Lemma}
\newtheorem{corollary}[theorem]{Corollary}
\newtheorem{proposition}[theorem]{Proposition}
\newtheorem{conjecture}{Conjecture}
\theoremstyle{definition}
\newtheorem{definition}[theorem]{Definition}
\newtheorem{question}[theorem]{Question}
\newtheorem{remark}[theorem]{Remark}
\newtheorem{remark*}[]{Remark}
\newtheorem{example}[theorem]{Example}
\newtheorem{example*}[]{Example}

\newtheorem{construction}[theorem]{Construction}
\newtheorem{construction*}[]{Construction}

\newtheorem{assumption}[theorem]{Assumption}
\newtheorem{assumption*}[]{Assumption}

\baselineskip=17pt
\parindent=0.5cm

\begin{abstract}

Extriangulated category was introduced by Nakaoka and Palu to give a unification of properties in exact categories and triangulated categories. A notion of tilting (or cotilting) subcategories in an extriangulated category is defined in this paper. We give a Bazzoni characterization of tilting (or cotilting) subcategories and obtain an Auslander-Reiten correspondence between tilting (cotilting) subcategories and coresolving covariantly (resolving contravariantly, resp.) finite subcatgories which are closed under direct summands and satisfies some cogenerating (generating, resp.) conditons. Applications of the results are given: we show that tilting (cotilting) subcategories defined here unify many previous works about tilting theory in module categories of Artin algebras and abelian categories admitting a cotorsion triples; we also show that the results work for triangulated categories with a proper class of triangles introduced by Beligiannis.  \\[0.5cm]
%\textbf{ 2010 Mathematics Subject Classification:}  18E10; 18E30; 18E40
\medskip
\end{abstract}

\section{Introduction}

Tilting modules or tilting functors (i.e. the functors induced by tilting modules as Hom-functors) as a generalization of Bernstein-Gelfand-Ponomarev reflection functors \cite{bgp} were introduced by Auslander, Platzeck and Reiten \cite{apr}, Brenner and Butler \cite{bb}, Happel and Ringel \cite{hr}. They have been played important roles and studied extensively in representation theory of algebras. In \cite{m}, Miyashita considered tilting modules of finite projective dimension, while Colby and the coauthors (\cite{ct}, \cite{cf}) studied generalized tilting modules of projective dimension one over arbitrary rings. Happel \cite{h} realized that tilting modules of finite projective dimension can be used to induce triangle equivalences between the algebras involved. Furthermore, Rickard \cite{r} gave a Morita theory for derived categories. For the case of tilting modules of finite projective dimension, an interesting correspondence between tilting (or cotilting) modules and covariantly (contravariantly, resp.) finite subcategories of the module category over an Artin algebra was gave in \cite{ar}. This result has important applications in the theory of algebraic groups via work on quasi-hereditary algebras \cite{r}. It has several generalizaions, see for example, \cite{ba}, \cite{dwzc}, \cite{as}.

The notion of extriangulated categories was introduced in \cite{np} as a simultaneous generalization of exact categories and triangulated categories. Exact categories and extension-closed subcategories of a triangulated category are extriangulated categories, while there are some other examples of extriangulated categories which are neither exact nor triangulated. Hence, many results hold on exact categories and triangulated categories can be unified in the same framework [\cite{inp}, \cite{ln}, \cite{np}, \cite{zz}]. Motivated by this idea, we introduce tilting subcategories in an extriangulated category and study their properties. This enables us to treat the tilting theory and its generalizations appeared before in a uniform way. More precisely, we obtain the Bazzoni's characterization of tilting (cotilting) subcategories and the Auslander-Reiten correspondence between tilting (cotilting) subcategories and coresolving covariantly (resolving contravariantly, resp.) finite subcatgories in the extriangulated category. As applications, our results unify that about tilting modules (or subcategories) of finite projective dimension. [\cite{ar}, \cite{ba}, \cite{dwzc}, \cite{as}]. We show that the relative tilting subcategories in an abelian category $\A$ with a complete cotorsion triple $(\X,\Y,\Z)$ studied in \cite{dwzc} are tilting subcategories in our sense. Therefore the results in \cite{dwzc} follow from our main results. We also show that a triangulated category $\C$ with a proper class of triangles $\mathcal{E}$ introduced by Beligiannis \cite{be} forms an extriangulated category. If in addition $(\C,\mathcal{E})$ has enough projectives and injectives, we deduce the same results on tilting (cotilting) subcategories in this setting.

The paper is organized as follows. In Section 2, we recall the definition of an extriangulated category and outline some basic properties that will be used later. In Section 3, we show how to construct covariantly finite or contravariantly finite subcategories from other contravariantly or covariantly finite subcategories in an extriangulated category. Then we give several subcategories associated to a self-orthogonal subcategory and illustrate some properties which are essential for our main results. In Section 4, we define tilting and cotilting subcategories in an extriangulated category, formulate the Bazzoni's characterization and establish a one-to-one correspondence between tilting (cotilting, resp.) subcategories and coresolving covariantly (resolving contravariantly, resp.) finite subcatgories which are closed under direct summands and satisfy some cogenerating (generating, resp.) conditions. This bijection is the Auslander-Reiten bijection established in \cite{ar} when the extriangulated category is the module category of finite generated left modules over an Artin algebra. In Section 5, we give two applications of our results: the first is the relative homology case established in \cite{dwzc}, the second is the triangulated categories with proper class of triangles defined by \cite{be}, we deduce the same results about tilting (cotilting) subcategories in these settings.

\section{Preliminaries}
Throughout the article, $\C$ denotes a Krull-Schmidt category, all subcategories are full additive subcategories closed under isomorphisms.
We denote by $\C(A,B)$ or $\Hom_{\C}(A,B)$ the set of morphisms from $A$ to $ B$ in $\C$.

We recall the definition of functorially finite subcategories in $\C$.
Let $\X$ be a subcategory of $\C$, a
morphism $f_B\colon X_B\to B$ of $\C$ with $X_B$ an object in $\X$ ($g_B\colon B\to Y_B$ of $\C$ with $Y_B$ an object in $\X$, resp.)
is said to be a {right $\X$-approximation} (left $\X$-approximation, resp.) of $B$, if the
morphisms $\C(X,f_B)\colon\C(X,X_B)\to\C(X,B)$ ($\C(g_B,X)\colon\C(B,X)\to\C(Y_B,X)$, resp.) are surjectives
for all objects $X$ in $\X$.
Moreover, one says that
$\X$ is contravariantly finite (covariantly finite, resp.) in $\C$ if every object in $\C$ has a right (left, resp.) $\X$-approximation, and functorially finite if it is both contravariantly finite and covariantly finite in $\C$.

We recall the definition and some basic properties of extriangulated categories from \cite{np}.

Let $\C$ be an additive category. Suppose that $\C$ is equipped with a biadditive functor $\E\colon\C\op\times\C\to\Ab$. For any pair of objects $A,C\in\C$, an element $\delta\in\E(C,A)$ is called an $\E$-extension. Thus formally, an $\E$-extension is a triplet $(A,\delta,C)$.
Let $(A,\del,C)$ be an $\E$-extension. Since $\E$ is a bifunctor, for any $a\in\C(A,A')$ and $c\in\C(C',C)$, we have $\E$-extensions
$$ \E(C,a)(\del)\in\E(C,A')\ \ \text{and}\ \ \ \E(c,A)(\del)\in\E(C',A). $$
We abbreviate denote them by $a_\ast\del$ and $c^\ast\del$ respectively.
For any $A,C\in\C$, the zero element $0\in\E(C,A)$ is called the spilt $\E$-extension.

\begin{definition}{[\cite{np}, Definition 2.3]}\label{mo}
Let $(A,\del,C),(A',\del',C')$ be any pair of $\E$-extensions. A morphism $$(a,c)\colon(A,\del,C)\to(A',\del',C')$$ of $\E$-extensions is a pair of morphisms $a\in\C(A,A')$ and $c\in\C(C,C')$ in $\C$, satisfying the equality
$$ a_\ast\del=c^\ast\del'. $$
Simply we denote it as $(a,c)\colon\del\to\del'$.
\end{definition}

Let $A,C\in\C$ be any pair of objects. Sequences of morphisms in $\C$
$$\xymatrix@C=0.7cm{A\ar[r]^{x} & B \ar[r]^{y} & C}\ \ \text{and}\ \ \ \xymatrix@C=0.7cm{A\ar[r]^{x'} & B' \ar[r]^{y'} & C}$$
are said to be equivalent if there exists an isomorphism $b\in\C(B,B')$ which makes the following diagram commutative.
$$\xymatrix{
A \ar[r]^x \ar@{=}[d] & B\ar[r]^y \ar[d]_{\simeq}^{b} & C\ar@{=}[d]&\\
A\ar[r]^{x'} & B' \ar[r]^{y'} & C &}$$

We denote the equivalence class of $\xymatrix@C=0.7cm{A\ar[r]^{x} & B \ar[r]^{y} & C}$ by $[\xymatrix@C=0.7cm{A\ar[r]^{x} & B \ar[r]^{y} & C}]$. For any $A,C\in\C$, we denote as
$ 0=[A\xrightarrow{\binom{1}{0}}A\oplus C\xrightarrow{(0,\ 1)}C].$

\begin{definition}{[\cite{np}, Definition 2.9]}\label{re}
For any $\E$-extension $\delta\in\E(C,A)$, one can associate an equivalence class $\s(\delta)=[\xymatrix@C=0.7cm{A\ar[r]^{x} & B \ar[r]^{y} & C}].$  This $\s$ is called a realization of $\E$, if it satisfies the following condition:
\begin{itemize}
\item Let $\del\in\E(C,A)$ and $\del'\in\E(C',A')$ be any pair of $\E$-extensions, with $$\s(\del)=[\xymatrix@C=0.7cm{A\ar[r]^{x} & B \ar[r]^{y} & C}],\ \ \ \s(\del')=[\xymatrix@C=0.7cm{A'\ar[r]^{x'} & B'\ar[r]^{y'} & C'}].$$
Then, for any morphism $(a,c)\colon\del\to\del'$, there exists $b\in\C(B,B')$ which makes the following diagram commutative.
\begin{equation}\label{reali}
\begin{array}{l}
$$\xymatrix{
A \ar[r]^x \ar[d]^a & B\ar[r]^y \ar[d]^{b} & C\ar[d]^c&\\
A'\ar[r]^{x'} & B' \ar[r]^{y'} & C' &}$$
\end{array}
\end{equation}
\end{itemize}
In this case, we say that sequence $\xymatrix@C=0.7cm{A\ar[r]^{x} & B \ar[r]^{y} & C}$ realizes \ $\del$, whenever it satisfies $\s(\del)=[\xymatrix@C=0.7cm{A\ar[r]^{x} & B \ar[r]^{y} & C}]$.
Remark that this condition does not depend on the choice of the representatives of the equivalence classes. In the above situation, we say that (\ref{reali}) (the triplet $(a,b,c)$, resp.) realizes $(a,c)$.
\end{definition}

\begin{definition}{[\cite{np}, Definition 2.10]}\label{rea}
Let $\C, \E$ be as above. A realization $\s$ of $\E$ is said to be additive if the following conditions are satisfied:
\begin{itemize}
\item[{\rm (i)}] For any $A,C\in\C$, the split $\E-$extension $0\in \E(C,A)$ satisfies $\s(0)=0$.
\item[{\rm (ii)}] For any pair of $\E-$extensions $\del=(A,\del,C)$ and $\del'=(A',\del',C')$,
                  $$\s(\del\oplus \del')=\s(\del)\oplus \s(\del')$$ holds.
\end{itemize}
\end{definition}

\begin{definition}{[\cite{np}, Definition 2.12]}\label{ext}
We call the pair $(\E,\s)$ an external triangulation of $\C$ if the following conditions are satisfied:
\begin{itemize}
\item[{\rm (ET1)}] $\E\colon\C\op\times\C\to\Ab$ is a biadditive functor.
\item[{\rm (ET2)}] $\s$ is an additive realization of $\E$.
\item[{\rm (ET3)}] Let $\del\in\E(C,A)$ and $\del'\in\E(C',A')$ be any pair of $\E$-extensions, realized as
$$ \s(\del)=[\xymatrix@C=0.7cm{A\ar[r]^{x} & B \ar[r]^{y} & C}],\ \ \s(\del')=[\xymatrix@C=0.7cm{A'\ar[r]^{x'} & B' \ar[r]^{y'} & C'}]. $$
For any commutative square
$$\xymatrix{
A \ar[r]^x \ar[d]^a & B\ar[r]^y \ar[d]^{b} & C&\\
A'\ar[r]^{x'} & B' \ar[r]^{y'} & C' &}$$
in $\C$, there exists a morphism $(a,c)\colon\del\to\del'$ which is realized by $(a,b,c)$.
\item[{\rm (ET3)$\op$}] Let $\del\in\E(C,A)$ and $\del'\in\E(C',A')$ be any pair of $\E$-extensions, realized by
$$\xymatrix@C=0.7cm{A\ar[r]^{x} & B \ar[r]^{y} & C}\ \ \text{and}\ \ \ \xymatrix@C=0.7cm{A'\ar[r]^{x'} & B' \ar[r]^{y'} & C'}$$
respectively.
For any commutative square
$$\xymatrix{
A \ar[r]^x& B\ar[r]^y \ar[d]^{b} & C\ar[d]^c&\\
A'\ar[r]^{x'} & B' \ar[r]^{y'} & C' &}$$
in $\C$, there exists a morphism $(a,c)\colon\del\to\del'$ which is realized by $(a,b,c)$.

\item[{\rm (ET4)}] Let $(A,\del,D)$ and $(B,\del',F)$ be $\E$-extensions realized by
$$\xymatrix@C=0.7cm{A\ar[r]^{f} & B \ar[r]^{f'} & D}\ \ \text{and}\ \ \ \xymatrix@C=0.7cm{B\ar[r]^{g} & C \ar[r]^{g'} & F}$$
respectively. Then there exist an object $E\in\C$, a commutative diagram
$$\xymatrix{A\ar[r]^{f}\ar@{=}[d]&B\ar[r]^{f'}\ar[d]^{g}&D\ar[d]^{d}\\
A\ar[r]^{h}&C\ar[d]^{g'}\ar[r]^{h'}&E\ar[d]^{e}\\
&F\ar@{=}[r]&F}$$
in $\C$, and an $\E$-extension $\del^{''}\in\E(E,A)$ realized by $\xymatrix@C=0.7cm{A\ar[r]^{h} & C \ar[r]^{h'} & E},$ which satisfy the following compatibilities.
\begin{itemize}
\item[{\rm (i)}] $\xymatrix@C=0.7cm{D\ar[r]^{d} & E \ar[r]^{e} & F}$  realizes $f'_{\ast}\del'$,
\item[{\rm (ii)}] $d^\ast\del''=\del$,

\item[{\rm (iii)}] $f_{\ast}\del''=e^{\ast}\del'$.
\end{itemize}

\item[{\rm (ET4)$\op$}]  Let $(D,\del,B)$ and $(F,\del',C)$ be $\E$-extensions realized by
$$\xymatrix@C=0.7cm{D\ar[r]^{f'} & A \ar[r]^{f} & B}\ \ \text{and}\ \ \ \xymatrix@C=0.7cm{F\ar[r]^{g'} & B \ar[r]^{g} & C}$$
respectively. Then there exist an object $E\in\C$, a commutative diagram
$$\xymatrix{D\ar[r]^{d}\ar@{=}[d]&E\ar[r]^{e}\ar[d]^{h'}&F\ar[d]^{g'}\\
D\ar[r]^{f'}&A\ar[d]^{h}\ar[r]^{f}&B\ar[d]^{g}\\
&C\ar@{=}[r]&C}$$
in $\C$, and an $\E$-extension $\del^{''}\in\E(C,E)$ realized by $\xymatrix@C=0.7cm{E\ar[r]^{h'} & A \ar[r]^{h} & C},$ which satisfy the following compatibilities.
\begin{itemize}
\item[{\rm (i)}] $\xymatrix@C=0.7cm{D\ar[r]^{d} & E \ar[r]^{e} & F}$  realizes $g'^{\ast}\del$,
\item[{\rm (ii)}] $\del'=e_\ast\del''$,

\item[{\rm (iii)}] $d_\ast\del=g^{\ast}\del''$.
\end{itemize}
\end{itemize}
In this case, we call $\s$ an $\E$-triangulation of $\C$, and call the triplet $(\C,\E,\s)$ an externally triangulated category, or for short,  extriangulated category.
\end{definition}

For an extriangulated category $\C$, we use the following notation:

\begin{itemize}
\item A sequence $\xymatrix@C=0.41cm{A\ar[r]^{x} & B \ar[r]^{y} & C}$ is called a conflation if it realizes some $\E$-extension $\del\in\E(C,A)$.
\item A morphism $f\in\C(A,B)$ is called an inflation if it admits some conflation $\xymatrix@C=0.7cm{A\ar[r]^{f} & B \ar[r]& C}.$
\item A morphism $f\in\C(A,B)$ is called a deflation if it admits some conflation $\xymatrix@C=0.7cm{K\ar[r]& A \ar[r]^f& B}.$
\item If a conflation $\xymatrix@C=0.6cm{A\ar[r]^{x} & B \ar[r]^{y} & C}$ realizes $\del\in\E(C,A)$, we call the pair $(\xymatrix@C=0.41cm{A\ar[r]^{x} & B \ar[r]^{y} & C},\del)$ an $\E$-triangle, and write it in the following way.
$$\xymatrix{A\ar[r]^{x} & B \ar[r]^{y} & C\ar@{-->}[r]^{\del}&}$$
\item Given an $\E-$triangle $\xymatrix{A\ar[r]^{x} & B \ar[r]^{y} & C\ar@{-->}[r]^{\del}&},$ we call $A$ the CoCone of $y\colon B\to C,$ and denote it by $CoCone(B\to C),$ also denote by $CoCone(y);$ we call $C$ the Cone of $x\colon A\to B,$ and denote it by $Cone(A\to B),$ also denote by $Cone(x).$
\item Let $\xymatrix{A\ar[r]^{x}&B\ar[r]^{y}&C\ar@{-->}[r]^{\delta}&}$ and $\xymatrix{A'\ar[r]^{x'}&B'\ar[r]^{y'}&C'\ar@{-->}[r]^{\delta'}&}$ be any pair of $\E$-triangles. If a triplet $(a,b,c)$ realizes $(a,c)\colon\del\to\del'$ as in $(\ref{t1})$, then we write it as
    $$\xymatrix{
A \ar[r]^x \ar[d]_{a} & B\ar[r]^y \ar[d]^{b} & C\ar@{-->}[r]^{\delta} \ar[d]^{c}&\\
A'\ar[r]^{x'} & B' \ar[r]^{y'} & C' \ar@{-->}[r]^{\delta'}&}$$
and call $(a,b,c)$ a morphism of $\E$-triangles.
\item A subcategory $\cal{T}$ of $\C$ is called extension-closed if $\cal{T}$ is closed under extensions, i.e. for any $\E$-triangle $\xymatrix{A\ar[r]^{x} & B \ar[r]^{y} & C\ar@{-->}[r]^{\del}&} $ with $A, C\in \cal{T}$, we have $B\in \cal{T}$.
\end{itemize}

\begin{remark}\label{sp}
For any objects $A,B\in\C$, we have  $\E$-triangles in $\C$: $$\xymatrix{A\ar[r]^{\binom{1}{0}\quad}&A\oplus B\ar[r]^{\quad(0,\ 1)}&B\ar@{-->}[r]^{0}&}\ \ \mbox{and} \ \ \xymatrix{A\ar[r]^{\binom{0}{1}\quad}&B\oplus A\ar[r]^{\quad(1,\, 0)}&B\ar@{-->}[r]^{0}&}.$$
\end{remark}

\begin{example}\label{exa1}
(1) Exact category and extension-closed subcategories of an triangulated category are extriangulated categories. The extension-closed subcategories of an extriangulated category are extriangulated categories, see \cite{np} for more detail.
\vspace{2mm}

(2) Let $\C$ be an extriangulated category, and $\cal J$ a subcategory of $\C$.
If $\cal J\subseteq\cal P\cap\cal I$, where $\cal P$ is the subcategory of projective objects in $\C$ and $\cal I$ is the subcategory of injective objects in $\C$ (see Definition \ref{pro}), then $\C/\cal J$ is an extriangulated category.
This construction gives extriangulated categories which are neither exact nor triangulated in general.
For more details, see [\cite{np}, Proposition 3.30]. There are also other such examples in \cite{zz}.
\end{example}

There are some basic results on extriangulated categories which are needed later on.

\begin{lemma}\label{lon}{\emph{[\cite{np}, Corollary 3.12]}}
Let $\C$ be an extriangulated category, $$\xymatrix{A\ar[r]^{x}&B\ar[r]^{y}&C\ar@{-->}[r]^{\delta}&}$$
an $\E$-triangle. Then we have the following long exact sequence:
$$\C(-, A)\xrightarrow{\C(-,x)}\C(-, B)\xrightarrow{\C(-,y)}\C(-, C)\xrightarrow{\delta^{\sharp}_-}
\E(-, A)\xrightarrow{\E(-,x)}\E(-, B)\xrightarrow{\E(-,y)}\E(-, C)
,$$
$$\C(C,-)\xrightarrow{\C(y,-)}\C(B,-)\xrightarrow{\C(x,-)}\C(A,-)\xrightarrow{\delta_{\sharp}^-}
\E(C,-)\xrightarrow{\E(y,-)}\E(B,-)\xrightarrow{\E(x,-)}\E(A,-)
.$$
\end{lemma}

\begin{lemma}\label{diag}{\emph{[\cite{np}, Corollary 3.15]}}
Let $(\C,\E,\s)$ be an extriangulated category. Then the following results hold.

\item{(1)}  Let $C$ be any object, and let
$$\xymatrix{A_{1}\ar[r]^{x_{1}}&B_{1}\ar[r]^{y_{1}}&C\ar@{-->}[r]^{\delta_{1}}&}, \xymatrix{A_{2}\ar[r]^{x_{2}}&B_{2}\ar[r]^{y_{2}}&C\ar@{-->}[r]^{\delta_{2}}&}$$
be any pair of $\E-$triangles. Then there is a commutative diagram in $\C$
$$\xymatrix{&A_{2}\ar@{=}[r]\ar[d]^{m_{2}}&A_{2}\ar[d]^{x_{2}}\\
A_{1}\ar[r]^{m_{1}}\ar@{=}[d]&M\ar[r]^{e_{1}}\ar[d]^{e_{2}}&B_{2}\ar[d]^{y_{2}}&\\
A_{1}\ar[r]^{x_{1}}&B_{1}\ar[r]^{y_{1}}&C&\\
&&}$$ which satisfies
$$\s(y_{2}^{*}\delta_{1})=[\xymatrix@C=0.41cm{A_{1}\ar[r]^{m_{1}} & M \ar[r]^{e_{1}} & B_{2}}],$$ $$\s(y_{1}^{*}\delta_{2})=[\xymatrix@C=0.41cm{A_{2}\ar[r]^{m_{2}}&M\ar[r]^{e_{2}}&B_{1}}],$$
$$m_{1*}\delta_{1}+m_{2*}\delta_{2}=0.$$
\item{(2)}  Dual of {(1)}.
\end{lemma}

The following lemma was given in [\cite{ln}, Proposition 1.20], which is another version of [NP, Corollary 3.16].
\begin{lemma}\label{pull}
Let $\xymatrix{A\ar[r]^{x}&B\ar[r]^{y}&C\ar@{-->}[r]^{\delta}&}$ be an $\E$-triangle, $f\colon A\to D$ be any morphism, and let $\xymatrix{D\ar[r]^{d}&E\ar[r]^{e}&C\ar@{-->}[r]^{f_{*}\delta}&}$ be any $\E$-triangle realizing $f_{*}\delta.$ Then there is a morphism $g$ which gives a morphism of $\E$-triangles
$$\xymatrix{A\ar[r]^{x}\ar[d]^f&B\ar[r]^y\ar[d]^g&C\ar@{-->}[r]^\del\ar@{=}[d]&\\
D\ar[r]^{d}&E\ar[r]^{e}&C\ar@{-->}[r]^{f_{*}\del}&}
$$
and moreover, $\xymatrix{A\ar[r]^{\binom{-f}{x}\quad}&D\oplus B\ar[r]^{\quad(d,\ g)}&E\ar@{-->}[r]^{e^{*}\delta}&}$ becomes an $\E$-triangle.
\end{lemma}

\begin{definition}{[\cite{np}, Definition 3.23]}\label{pro}
Let $\C, \E$ be as above. An object $P\in \C$ is called projective if it satisfies the following condition.
\begin{itemize}
\item For any $\E-$triangle $\xymatrix@C=0.41cm{A\ar[r]^{x} & B \ar[r]^{y}&C\ar@{-->}[r]^{\delta}&}$ and any morphism $c\in \C(P,C)$, there exists $b\in \C(P,B)$ satisfying $y\circ b=c$.
\end{itemize}
Injective objects are defined dually.

We denote the subcategory consisting of projective objects in $\C$ by $Proj(\C)$.
Dually, the subcategory of injective objects in $\C$ is denoted by $Inj(\C)$.
\end{definition}

\begin{lemma}{[\cite{np}, Proposition 3.24]}\label{proj}
An object $P\in \C$ is projective if and only if it satisfies $\E(P,A)=0$ for any $A\in \C$.
The dual property also holds for injective objects.
\end{lemma}

\begin{definition}{[\cite{np}, Definition 3.25]}\label{proje}
Let $(\C,\E,\s)$ be an extriangulated category. We say that it has enough projectives (enough injective, resp.) if it satisfies the following condition:
\begin{itemize}
\item For any object $C\in \C$ ($A\in \C$, resp.), there exists an $\E-$triangle $$\xymatrix{A\ar[r]^{x}&P\ar[r]^{y}&C\ar@{-->}[r]^{\delta}&} (\xymatrix{A\ar[r]^{x}&I\ar[r]^{y}&C\ar@{-->}[r]^{\delta}&}, \mbox{resp.})$$ satisfying $P\in Proj(\C)$ ($I\in Inj(\C)$, resp.).

In this case, $A$ is called the syzygy of $C$ ($C$ is called the cosyzygy of $A$, resp.) and is denoted by $\Omega(C) (\Sigma(A),$ resp.).
\end{itemize}
\end{definition}

Suppose $\C$ is an extriangulated category with enough projectives and injectives. For a subcategory $\B\subseteq\C$, put $\Omega^{0}\B=\B$, and for $i>0$ we define $\Omega^{i}\B$ inductively to be the subcategory consisting of syzygies of objects in $\Omega^{i-1}\B$, i.e.
$$\Omega^{i}\B=\Omega(\Omega^{i-1}\B).$$
We call $\Omega^{i}\B$ the $i$-th syzygy of $\B$. Dually we define the $i$-th cosyzygy $\Sigma^{i}\B$ by $\Sigma^{0}\B=\B$ and $\Sigma^{i}\B=\Sigma(\Sigma^{i-1}\B)$ for $i>0$.

In [\cite{ln} 5.2] the authors defined higher extension groups in an extriangulated category having enough projectives and injectives as $\E^{i+1}(X,Y)\cong\E(X,\Sigma^{i}Y)\cong \E(\Omega^{i}X,Y)$ for $i\geq0,$ and they showed the following result:

\begin{lemma}{[\cite{ln}, Proposition 5.2]}\label{long}
Let $\xymatrix{A\ar[r]^{x}&B\ar[r]^{y}&C\ar@{-->}[r]^{\delta}&}$ be an $\E-triangle$. For any object $X\in \B$, there are long exact sequences
$$\cdots\rightarrow\E^{i}(X, A)\xrightarrow{x_{*}}\E^{i}(X, B)\xrightarrow{y_{*}}\E^{i}(X, C)\rightarrow\E^{i+1}(X, A)\xrightarrow{x_{*}}\E^{i+1}(X, B)\xrightarrow{y_{*}}\cdots (i\geq 1),$$
$$\cdots\rightarrow\E^{i}(C, X)\xrightarrow{y^{*}}\E^{i}(B, X)\xrightarrow{x^{*}}\E^{i}(A, X)\rightarrow\E^{i+1}(C, X)\xrightarrow{y^{*}}\E^{i+1}(B, X)\xrightarrow{x^{*}}\cdots (i\geq 1).$$
\end{lemma}

An $\E-$triangle sequence in $\C$ is displayed as a sequence
$$\cdots\rightarrow X_{n+1}\xrightarrow{d_{n+1}}X_{n}\xrightarrow{d_{n}}X_{n-1}\rightarrow\cdots$$
over $\C$ such that for any $n$, there are $\E-$triangles $\xymatrix{K_{n+1}\ar[r]^{g_{n}}&X_{n}\ar[r]^{f_{n}}&K_{n}\ar@{-->}[r]^{\delta^{n}}&}$ and the differential $d_{n}=g_{n-1}f_{n}.$

Let $\C$ be an extriangulated category with enough projectives and injectives. Then for any $A\in\C$, we have an $\E-$triangle $\xymatrix{A_{1}\ar[r]^{g_{0}}&P_{0}\ar[r]^{f_{0}}&A\ar@{-->}[r]^{\delta_{0}}&}$ with $P_{0}\in Proj(\C).$
For $A_{1},$ we have an $\E-$triangle $\xymatrix{A_{2}\ar[r]^{g_{1}}&P_{1}\ar[r]^{f_{1}}&A_{1}\ar@{-->}[r]^{\delta_{1}}&}$ with $P_{1}\in Proj(\C).$ Inductively, one can see that $A$ admits a projective resolution as follows:
$$\cdots\rightarrow P_{n}\xrightarrow{d_{n}}P_{n-1}\rightarrow\cdots  \rightarrow P_{1}\xrightarrow{d_{1}}P_{0}\xrightarrow{d_{0}}A $$
%$$A\xrightarrow{d^{0}}I^{0}\xrightarrow{d^{1}}I^{1}\rightarrow\cdots  \rightarrow I^{n-1}\xrightarrow{d^{n}}I^{n}\rightarrow\cdots $$
with $P_{i}\in Proj(\C),$ and $d_{n}=g_{n-1}f_{n}.$
We call a projective resolution of $A$ is of length $n$ if $A_n=P_n$ is projective and $d_n=g_{n-1}$ is an inflation. We define the projective dimension of $A$, denoted by $pd(A)$, the minimal length $n$ of all projective resolutions of $A$, if there is a projective resolution of $A$ of finite length; otherewise, we define the projective dimension of $A$ to be $\infty$. Dually, we define the injective resolution and injective dimension of $A\in\C,$ and the injective dimension of $A$ is denoted by $id(A).$

\begin{lemma}\label{projec}
Let $\C$ be an extriangulated category with enough projectives and injectives. Then for any $A\in\C,$ the following statements are equivalent:
\emph{(1)} $\emph{pd}A\leq{n};$
\emph{(2)} $\E^{n+1}(A,X)=0,$ for any  $X\in\C;$
\emph{(3)} $\E^{n+i}(A,X)=0,$ for any  $X\in\C$ and $i\geq1.$
\end{lemma}

\proof  $(3)\Rightarrow(2)$ is obvious.

$(1)\Rightarrow(2),(3).$ Let $$ P_{n}\xrightarrow{d_{n}}P_{n-1}\rightarrow\cdots  \rightarrow P_{1}\xrightarrow{d_{1}}P_{0}\xrightarrow{d_{0}}A $$ be a projective resolution of $A,$ then we have $\E^{n+i}(A,X)=\E^{i}(P_{n},X)=0,$ for any $X\in\C$ and $i\geq1.$

$(2)\Rightarrow(1).$ For any $A\in \C,$ let $$\cdots\rightarrow P_{n}\xrightarrow{d_{n}}P_{n-1}\rightarrow\cdots  \rightarrow P_{1}\xrightarrow{d_{1}}P_{0}\xrightarrow{d_{0}}A$$ be a projective resolution of $A.$ It follows that $\E^{n+i}(A,X)=\E^{i}(A_{n},X)=0,$ for any $X\in\C,$ where $A_{n}=\Omega^{n}A$. Applying Lemma \ref{proj} one can see $A_{n}\in Proj(\C),$ then $(1)$ holds.\qed

\begin{remark}\label{inj}
The dual result also holds for $\emph{id}A.$
\end{remark}

\section{Basic results}
\setcounter{equation}{0}

In this section we firstly show how to construct covariantly finite or contravariantly finite subcategories from other contravariantly or covariantly finite subcategories, then we outline some basic properties of self-orthonogal subcatgories that are essential for the following sections. From now on, we suppose all extriangulated categories have enough projectives and injectives.

\subsection{Covariantly and contravariantly finite subcategories}
For a subcategory $\X\subset\C,$ define $\X^{\bot_{1}}=\{Y\in\C\ | \ \E(X,Y)=0, \forall X\in\X \},$ \ $\X^{\bot}=\{Y\in\C \ | \ \E^{i}(X,Y)=0 \ , \forall i\geq1, \ X\in\X\}.$ Similarly, define $\sideset{^{\bot_{1}}}{}{\mathop{\X}}=\{Y\ \in\C\ |\ \E(Y,X)=0, \forall X\in\X \},$ \ $\sideset{^{\bot}}{}{\mathop{\X}}=\{Y\in\C\ |\ \E^{i}(Y,X)\ =0, \forall i\geq1, \ X\in\X\}.$

\begin{lemma}\label{cov}
Let $\X \subseteq\C$ be a covariantly finite subcategory closed under extensions and direct summands, and $\Y=\sideset{^{\bot_{1}}}{}{\mathop{\X}}.$ Then $\Y$ is a contravariantly finite subcategory closed under extensions and direct summands.

\proof  From the definition, it is obvious that $\Y$ is closed under extensions and direct summands.

 Now to show $\Y$ is contravariantly finite. For any $C\in\C,$ there is an $\E-$triangle $\xymatrix{C^{'}\ar[r]^{x}&P\ar[r]^{y}&C\ar@{-->}[r]^{\delta}&}$ with $P\in Proj(\C)$ since $\C$ has enough projectives. Applying Lemma \ref{lon} and Lemma \ref{proj}, we have the following exact sequence of functors: $$\C(C,-)\xrightarrow{\C(y,-)}\C(P,-)\xrightarrow{\C(x,-)}\C(C^{'},-)\xrightarrow{\delta_{\sharp}^-}
\E(C,-)\rightarrow0.$$
Then restricting to $\X,$ we have $$\C(C,-)|_{\X}\xrightarrow{\C(y,-)|_{\X}}\C(P,-)|_{\X}\xrightarrow{\C(x,-)|_{\X}}\C(C^{'},-)|_{\X}\xrightarrow{\delta_{\sharp}^-|_{\X}}
\E(C,-)|_{\X}\rightarrow0.$$
Since $\X$ is a covariantly finite subcategory, $C^{'}$ admits a left $\X-$approximation $f\colon C^{'}\to X^{'},$ we have a exact sequence $\C(X^{'},-)|_{\X}\xrightarrow{\C(f,-)|_{\X}}\C(C^{'},-)|_{\X}\rightarrow0.$
Composting $\C(f,-)|_{\X}$ with $\delta_{\sharp}^-|_{\X},$ we have the following exact sequence of functors: $$\C(X^{'},-)|_{\X}\rightarrow \E(C,-)|_{\X}\rightarrow0.$$
Since $\C$ is Krull-Schmit, we can take a projective cover $\phi\colon \C(X,-)|_{\X}\to \E(C,-)|_{\X}\rightarrow0$ of $\X-$modules (for the definition of $\X-$modules see \cite{a}).
This is induced by an $\E-$triangle $\xymatrix{X\ar[r]^{a}&M\ar[r]^{b}&C\ar@{-->}[r]^{\theta}&}.$
Suppose now we have shown that $M\in\Y,$ one can see the sequence $\C(\Y,M)\rightarrow\C(\Y,C)\rightarrow\E(\Y,X)=0$ is exact, then $C$ admits a right $\Y-$approximation.
So we now turn our attention to show that $M\in\Y,$ hence to show that any $\E-$triangle $\xymatrix{X_{1}\ar[r]^{c}&N\ar[r]^{d}&M\ar@{-->}[r]^{\delta_{1}}&}$ with $X_{1}\in\X$ splits.
Using {\rm (ET4)$\op$} we obtain the following commutative diagram:
$$\xymatrix{X_{1}\ar[r]^{m}\ar@{=}[d]&X_{2}\ar[r]^{n}\ar[d]^{g}&X\ar[d]^{a}\\
X_{1}\ar[r]^{c}&N\ar[d]^{h}\ar[r]^{d}&M\ar[d]^{b}\\
&C\ar@{=}[r]&C}$$
with $\s (a^{*}\delta_{1})=[\xymatrix{X_{1}\ar[r]^{m}&X_{2}\ar[r]^{n}&X}].$
Therefore, $X_{2}\in\X$ since $X,X_{1}\in \X$ and $\X$ is closed under extensions.
We have the commutative diagram
 $$\xymatrix{\C(C,-)|_{\X}\ar[r]^{\C(b,-)|_{\X}}\ar@{=}[d]&\C(M,-)|_{\X}\ar[r]^{\C(a,-)|_{\X}}\ar[d]^{\C(d,-)|_{\X}}&\C(X,-)|_{\X}\ar[r]^{\phi\quad}\ar[d]^{\C(n,-)|_{\X}}& \E(C,-)|_{\X}\rightarrow0\ar@{=}[d]\\
\C(C,-)|_{\X}\ar[r]^{\C(h,-)|_{\X}}&\C(N,-)|_{\X}\ar[r]^{\C(g,-)|_{\X}}&\C(X_{2},-)|_{\X}\ar[r]^{\psi\quad}&
\E(C,-)|_{\X}}$$
of exact sequences of $\X-$modules. As $\phi$ is a projective cover, $\C(n,-)|_{\X}$ is a split monomorphism, hence $n$ is a split monomorphism. Evaluating the upper sequence above at $X_{1},$ one can see the sequence $\xymatrix{0\rightarrow\E(M,X_{1})\ar[r]^{a^{*}}&\E(X,X_{1})}$ is exact. Therefore $\delta_{1}=0$ since $a^{*}\delta_{1}=[\xymatrix{X_{1}\ar[r]^{m}&X_{2}\ar[r]^{n}&X}]=0.$ Hence $\E(M,X_{1})=0$ for any $X_{1}\in\X,$ as desired. \qed
\end{lemma}
We state the dual result of Lemma \ref{cov} without proof.

\begin{lemma}\label{con}
Let $\Y \subseteq\C$ be a contravariantly finite subcategory closed under extensions and direct summands, and $\X=\Y^{\bot_{1}}.$ Then $\X$ is a covariantly finite subcategory closed under extensions and direct summands.
\end{lemma}

From the proof of Lemma \ref{cov} and the dual version, we have the following results.
\begin{corollary}\label{cocon}
\begin{itemize}
\item[(1)] Let $\X \subseteq\C$ be a covariantly finite subcategory closed under extensions and direct summands. Denote $\sideset{^{\bot_{1}}}{}{\mathop{\X}}$ by $\Y.$ Then for any $C\in\C,$ there is an $\E-$triangle $\xymatrix{X\ar[r]^{x}&Y\ar[r]^{y}&C\ar@{-->}[r]^{\delta}&}$ such that $X\in\X, Y\in\Y$ and $y\colon Y\to C$ is a right $\Y-$approximation of $C.$
\item[(2)] Let $\Y \subseteq\C$ be a cotravariantly finite subcategory closed under extensions and direct summands. Denote $\Y^{\bot_{1}}$ by $\X$. Then for any $C\in\C,$ there is an $\E-$triangle $\xymatrix{C\ar[r]^{x}&X\ar[r]^{y}&Y\ar@{-->}[r]^{\delta}&}$ such that $X\in\X, Y\in\Y$ and $x\colon C\to X$ is a left $\X-$approximation of $C.$
\end{itemize}
\end{corollary}

\subsection{Self-orthonogal subcategories}
The subcategory $\mathcal T$ is said to be self-orthonogal provided that $\E^{i}(T_{1},T_{2})=0, \forall i\geq1,$ and for any $T_{1},T_{2}\in\mathcal T.$
The symbol $\hat{\mathcal T_{n}}$ (or $\check{\mathcal T_{n}})$ denotes the subcategory of objects $A\in\C$ such that there exists an $\E-$triangle sequence $T_{n}\rightarrow T_{n-1}\cdots\rightarrow T_{0}\rightarrow A$ ($A\rightarrow T_{0}\cdots\rightarrow T_{n-1}\rightarrow T_{n}$, \mbox{resp.}) with each $T_{i}\in\mathcal T.$
We denote by $\hat{\mathcal T}$(or $\check{\mathcal T})$ the union of all $\hat{\mathcal T_{n}}$ $(\check{\mathcal T_{n}}, \mbox{resp.})$ for some non-negative $n.$ That is to say
$$\hat{\mathcal T}=\bigcup_{n=0}^{\infty}\hat{\mathcal T_{n}},\ \check{\mathcal T}=\bigcup_{n=0}^{\infty}\check{\mathcal T_{n}}.$$

\begin{lemma}\label{ort}
Suppose $\mathcal T$ is a self-orthonogal subcategory of $\C.$ Then $\E^{i}(X_{1},X_{2})=0$, for any $i\geq1$ and for any $X_{1}\in\check{\mathcal T}$, $X_{2}\in\mathcal T^{\bot}.$
\end{lemma}

\proof Since $X_{1}\in\check{\mathcal T},$ there is a $n,$ such that $X_{1}\in\check{\mathcal T_{n}}.$ We obtain an $\E-$triangle sequence
$$X_{1}\rightarrow T^{0}\cdots\rightarrow T^{n-1}\rightarrow T^{n}$$ with each $T^{i}\in\mathcal T.$ Therefore, we have $\E^{i}(X_{1},X_{2})=\E^{i+n}(T^{n},X_{2})=0, \forall i\geq1.$\qed

Given a self-orthonogal subcategory $\mathcal T.$
Note that $\mathcal T^{\bot}$ ($\sideset{^{\bot}}{}{\mathop{\mathcal T}}$, resp.) is the largest subcategory of $\C$ such that $\E^{i}(\mathcal T,\mathcal T^{\bot})=0, \forall i\geq1$ ($\E^{i}(\sideset{^{\bot}}{}{\mathop{\mathcal T}},\mathcal T)=0, \forall i\geq1$, resp.).
Instead of considering $\mathcal T^{\bot}$ (or $\sideset{^{\bot}}{}{\mathop{\mathcal T}}$), we introduce a subcategory ${}_\mathcal T\X \subset\mathcal T^{\bot}$ ($\X_{\mathcal T}\subset\sideset{^{\bot}}{}{\mathop{\mathcal T}}$, resp.) consists of objects $A\in\C$ which are generated (cogenerated, resp.) by objects of $\mathcal T$ in the sense that there is an $\E-$triangle sequence $\cdots\rightarrow T_{n}\xrightarrow{d_{n}}T_{n-1}\rightarrow\cdots  \rightarrow T_{1}\xrightarrow{d_{1}}T_{0}\xrightarrow{d_{0}}A$ with $T_{i}\in\mathcal T$ and $CoCone(d_{i})\in\mathcal T^{\bot}$ ($A\xrightarrow{d_{0}} T_{0}\xrightarrow{d_{1}}T_{1}\rightarrow\cdots  \rightarrow T_{n-1}\xrightarrow{d_{n}}T_{n}\rightarrow\cdots$ with $T_{i}\in\mathcal T$ and $Cone(d_{i})\in\sideset{^{\bot}}{}{\mathop{\mathcal T}}$, resp.).
It is obvious that $\mathcal T\subset {}_\mathcal T\X$ ($\mathcal T\subset \X_{\mathcal T}$) and ${}_\mathcal T\X$ ($\X_{\mathcal T}$, resp.) is the largest subcategory of $\C$ such that $\mathcal T$ is projective (injective, resp.) and a generator (cogenerator, resp.) in it.
It is clearly that $$\hat{\mathcal T}\subset {}_\mathcal T\X\subset\mathcal T^{\bot}, \ \check{\mathcal T}\subset\X_{\mathcal T}\subset\sideset{^{\bot}}{}{\mathop{\mathcal T}}.$$

Now we give a characterization of objects in $\check{\mathcal T_{n}}.$ The dual result also holds for $\hat{\mathcal T_{n}}.$

\begin{lemma}\label{orth}
Let $\mathcal T\in\C$ be a self-orthonogal subcategory closed under direct summands. For any object $X\in\X_{\mathcal T},$ the following statements are equivalent:
\begin{itemize}
\item[$(1)$] $X\in\check{\mathcal T_{n}};$

\item[$(2)$] $E^{n+1}(Y,X)=0, \forall \ Y\in\sideset{^{\bot}}{}{\mathop{\mathcal T}};$

\item[$(3)$] $E^{n+1}(Y,X)=0, \forall \ Y\in\X_{\mathcal T}.$
\end{itemize}
\end{lemma}

\proof The proof is the same as that of Lemma \ref{projec}, we left it to the reader.\qed

\begin{lemma}\label{ortho}
For a self-orthonogal subcategory $\mathcal T,$ ${}_\mathcal T\X$ is closed under extensions, direct summands and Cone of inflations.
\end{lemma}

\proof We first show that ${}_\mathcal T\X$ is closed under extensions.
Let $\xymatrix{A\ar[r]^{x}&B\ar[r]^{y}&C\ar@{-->}[r]^{\delta}&}$ be an $\E-$triangle with $A, C \in {}_\mathcal T\X.$ Then we have $\E-$triangles $$\xymatrix{K_{0}\ar[r]^{a_{0}}&T_{0}\ar[r]^{b_{0}}&A\ar@{-->}[r]^{\delta_{0}}&}$$ $$\xymatrix{K_{1}\ar[r]^{a_{1}}&T_{1}\ar[r]^{b_{1}}&C\ar@{-->}[r]^{\delta_{1}}&}$$
with $T_{0}, \ T_{1}\in \mathcal T$ and $K_{0}, \ K_{1}\in {}_\mathcal T\X.$ Using Lemma \ref{diag}, we have the following commutative diagram
$$\xymatrix{&A\ar@{=}[r]\ar[d]^{f}&A\ar[d]^{x}\\
K_{1}\ar[r]^{m}\ar@{=}[d]&U\ar[r]^{n}\ar[d]^{g}&B\ar[d]^{y}&\\
K_{1}\ar[r]^{a_{1}}&T_{1}\ar[r]^{b_{1}}&C&\\
&&}$$
We conclude that $\E(T_{1},A)=0$ since $T_{1}\in\mathcal T$ and $A\in {}_\mathcal T\X \subset\mathcal T^{\bot}.$ Then $U\cong A\oplus T_{1}.$
Using {\rm (ET4)$\op$}, we have the following commutative diagram
$$\xymatrix{K_{0}\ar[r]^{a}\ar@{=}[d]&V\ar[r]^{b}\ar[d]^{c}&K_{1}\ar[d]^{m}\\
K_{0}\ar[r]^{\binom{a_{0}}{0}}&T_{0}\oplus T_{1}\ar[d]^{d}\ar[r]^{\left(\begin{array}{cc}b_{0}&0\\
0&1\end{array}\right)}&A\oplus T_{1}\ar[d]^{n}\\
&B\ar@{=}[r]&B}$$
It follows that $V\in \mathcal T^{\bot}$ since $K_{0}, K_{1}\in\mathcal T^{\bot}$ and $\mathcal T^{\bot}$ is closed under extensions. Since $K_{0}, K_{1}\in {}_\mathcal T\X,$ one can consider $V$ instead of $B$ and repeat this process to show that $B\in {}_\mathcal T\X.$

Secondly, to show ${}_\mathcal T\X$ is closed under Cone of inflations.
Let $\xymatrix{A\ar[r]^{x}&B\ar[r]^{y}&C\ar@{-->}[r]^{\delta}&}$ be an $\E-$triangle with $A, B \in {}_\mathcal T\X.$ Then we have the following commutative diagram:
$$\xymatrix{K_{0}\ar[r]^{a}\ar@{=}[d]&E\ar[r]^{b}\ar[d]^{c}&A\ar[d]^{x}\\
K_{0}\ar[r]^{a_{0}}&T_{0}\ar[d]^{d}\ar[r]^{b_{0}}&B\ar[d]^{y}\\
&C\ar@{=}[r]&C}$$
with $T_{0}\in\mathcal T$ and $K_{0}\in {}_\mathcal T\X.$ One has $E\in {}_\mathcal T\X$ since $A\in {}_\mathcal T\X$ and ${}_\mathcal T\X$ is closed under extensions. Hence, repeating this process we have $C\in {}_\mathcal T\X.$

Now to show ${}_\mathcal T\X$ is closed under summands. Let $B=A\oplus C$ such that $B\in {}_\mathcal T\X.$ Then we have an $\E-$triangle $\xymatrix{K_{0}\ar[r]^{a_{0}}&T_{0}\ar[r]^{b_{0}}&B\ar@{-->}[r]^{\delta}&}$ with $T_{0}\in\mathcal T, K_{0}\in {}_\mathcal T\X.$ Using {\rm (ET4)$\op$}, we have the following commutative diagram
$$\xymatrix{K_{0}\ar[r]^{a}\ar@{=}[d]&U\ar[r]^{b}\ar[d]^{c}&A\ar[d]^{\binom{1}{0}}\\
K_{0}\ar[r]^{a_{0}}&T_{0}\ar[d]^{d}\ar[r]^{b_{0}}&B\ar[d]^{(0,\, 1)}\\
&C\ar@{=}[r]&C}$$
Then we have an $\E-$triangle $\xymatrix{K_{0}\ar[r]^{\binom{a}{0}}&U\oplus C\ar[r]^{\left(\begin{array}{cc}b&0\\
0&1\end{array}\right)}&A\oplus C\ar@{-->}[r]^{\theta}&}$ with $U\oplus C\in {}_\mathcal T\X$ since $K_{0}, A\oplus C\in {}_\mathcal T\X$ and ${}_\mathcal T\X$ is closed under extensions.
Hence $U$ is a direct summand of an object in ${}_\mathcal T\X.$ By repeating this process, we can show $C\in {}_\mathcal T\X,$ so ${}_\mathcal T\X$ is closed under direct summands.\qed

In general, ${}_\mathcal T\X\subsetneq\mathcal T^{\bot}.$ In what follows we give a characterization of when they are equal.
We call a subcategory $\X\subset\C$ a generator (cogenerator, resp.) for a subcategory $\Y\subset\C$ if $\X\subset\Y$ and for any object $Y\in\Y,$ there is an $\E-$triangle $\xymatrix{Y'\ar[r]^{x}&X\ar[r]^{y}&Y\ar@{-->}[r]^{\delta}&}$ ( $\xymatrix{Y\ar[r]^{x}&X\ar[r]^{y}&Y'\ar@{-->}[r]^{\delta}&}$, resp.) in $\C$ with $X\in\X$ and $Y'\in\Y.$

\begin{lemma}\label{orthon}
For a self-orthonogal subcategory $\mathcal T.$ The following statements are equivalent:
\begin{itemize}
\item[$(1)$] $\mathcal T$ is a generator for $\mathcal T^{\bot};$

\item[$(2)$] ${}_\mathcal T\X=\mathcal T^{\bot}.$
\end{itemize}
\end{lemma}

\proof $(1)\Rightarrow(2).$ We just need to show $\mathcal T^{\bot}\subset{}_\mathcal T\X.$ Since $\mathcal T$ is a generator for $\mathcal T^{\bot},$ for any $X\in \mathcal T^{\bot},$ there exists an $\E-$triangle $\xymatrix{X'\ar[r]^{x}&T\ar[r]^{y}&X\ar@{-->}[r]^{\delta}&}$ with $T\in\mathcal T , X'\in\mathcal T^{\bot}.$
Repeating this process by replacing $X'$ with $X,$ one can see $X\in {}_\mathcal T\X.$

$(2)\Rightarrow(1).$ Obvious from the definition.\qed

The following two lemmas are extriangulated category version of [DWZC lemma3.3, corollary3.4], we give a proof here for the convenience of the reader.
\begin{lemma}\label{ch}
Let $\X\subset\Y$ be subcategories of $\C$ such that $\Y$ is closed under extensions. Assume that $\X$ is a generator of $\Y.$ If there is an $\E-$triangle sequence
$$A'\rightarrow Y_{n}\rightarrow Y_{n-1}\cdots\rightarrow Y_{1}\rightarrow A$$
in $\C$ with each $Y_{i}\in\Y.$ Then there exist objects $U_{n}, V_{n}\in \C$ satisfying the following conditions:
\begin{itemize}
\item[$(1)$] $U_{n}\in\Y;$

\item[$(2)$] there exist $\E-$triangle sequences $\xymatrix{U_{n}\ar[r]^{x}&V_{n}\ar[r]^{y}&A'\ar@{-->}[r]^{\delta}&}$ and
$$V_{n}\rightarrow X_{n}\rightarrow X_{n-1}\cdots\rightarrow X_{1}\rightarrow A$$
in $\C$ with each $X_{i}\in\X.$
\end{itemize}
\end{lemma}

\proof We proceed by induction on $n.$

If $n=1,$ then we have an $\E-$triangle $\xymatrix{A'\ar[r]^{x'}&Y_{1}\ar[r]^{y'}&A\ar@{-->}[r]^{\delta_{1}}&}$ in $\C$ with $Y_{1}\in \Y.$
Since $\X$ is a generator of $\Y,$ there is an $\E-$triangle $\xymatrix{U_{1}\ar[r]^{x_{1}}&X_{1}\ar[r]^{y_{1}}&Y_{1}\ar@{-->}[r]^{\delta_{2}}&}$ in $\C$ with $X_{1}\in \X$ and $U_{1}\in\Y.$ Then we have the following commutative diagram:
$$\xymatrix{U_{1}\ar[r]\ar@{=}[d]&V_{1}\ar[r]\ar[d]&A'\ar[d]\\
U_{1}\ar[r]&X_{1}\ar[d]\ar[r]&Y_{1}\ar[d]\\
&A\ar@{=}[r]&A}$$
the first row and the second column are the desired $\E-$triangles.

Now suppose that the result hold for $n-1.$ Let $A''=Cone(A'\to Y_{n}),$ then by the inductive assumption, there are $\E-$triangles:
$$\xymatrix{U_{n-1}^{'}\ar[r]^{x'}&V_{n-1}^{'}\ar[r]^{y'}&A''\ar@{-->}[r]^{\delta'}&}$$
and
$$V_{n-1}^{'}\rightarrow X_{n-1}\rightarrow X_{n-2}\cdots\rightarrow X_{1}\rightarrow A$$
in $\C$ with $U_{n-1}^{'}\in\Y$ and $X_{i}\in\X.$
Then using Lemma \ref{diag}, we have the following commutative diagram:
$$\xymatrix{&U_{n-1}^{'}\ar@{=}[r]\ar[d]&U_{n-1}^{'}\ar[d]\\
A'\ar[r]\ar@{=}[d]&M\ar[r]\ar[d]&V_{n-1}^{'}\ar[d]&\\
A'\ar[r]&Y_{n}\ar[r]&A''&\\
&&}$$ One can see that $M\in\Y$ since $\Y$ is closed under extensions and $U_{n-1}^{'}, Y_{n}\in \Y$, hence there is an $\E-$triangle $\xymatrix{U_{n}\ar[r]^{a}&X\ar[r]^{b}&M\ar@{-->}[r]^{\theta}&}$ with $X\in\X$ and $U_{n}\in\Y.$ Therefore, there is a commutative diagram:
$$\xymatrix{U_{n}\ar[r]\ar@{=}[d]&V_{n}\ar[r]\ar[d]&A'\ar[d]\\
U_{n}\ar[r]&X\ar[d]\ar[r]&M\ar[d]\\
&V_{n-1}^{'}\ar@{=}[r]&V_{n-1}^{'}}$$
One has that the first row is the desired $\E-$triangle and $U_{n}, V_{n}$ are the desired objects.\qed

\begin{lemma}\label{chang}
Let $\X\subset\Y$ be subcategories of $\C$ such that $\Y$ is closed under extensions. Assume that $\X$ is a generator of $\Y.$ If $A'\in\check{\Y_{n}},$ then there exists an $\E-$triangle
$$\xymatrix{U\ar[r]^{x}&V\ar[r]^{y}&A'\ar@{-->}[r]^{\delta}&}$$
in $\C$ with $U\in\Y$ and $V\in\check{\X_{n}}.$
\end{lemma}

\proof Since $A'\in\check{\Y_{n}},$ we have the following $\E-$triangle sequence:
$$A'\rightarrow Y_{0}\rightarrow Y_{1}\cdots\rightarrow Y_{n-1}\rightarrow Y_{n}$$
in $\C$ with $Y_{i}\in\Y.$
Since $\X$ is a generator of $\Y,$ there is an $\E-$triangle $\xymatrix{Y_{n}^{'}\ar[r]^{x'}&X_{0}\ar[r]^{y'}&Y_{n}\ar@{-->}[r]^{\delta}&}$ with $X_{0}\in\X$ and $Y_{n}^{'}\in\Y.$
Let $A_{1}=CoCone(Y_{n-1}\to Y_{n}),$ then using Lemma \ref{diag} we have the following commutative diagram:
$$\xymatrix{&A_{1}\ar@{=}[r]\ar[d]&A_{1}\ar[d]\\
Y_{n}^{'}\ar[r]\ar@{=}[d]&E\ar[r]\ar[d]&Y_{n-1}\ar[d]&\\
Y_{n}^{'}\ar[r]&X_{0}\ar[r]&Y_{n}&\\&&}$$
Hence we have the following $\E-$triangle sequence:
$$A'\rightarrow Y_{0}\rightarrow Y_{1}\cdots\rightarrow Y_{n-2}\rightarrow E\rightarrow X_{0}$$
in $\C.$
One can see that $E\in\Y$ since $Y_{n}^{'}, Y_{n-1}\in\Y$ and $\Y$ is closed under extensions.
Hence, using Lemma \ref{ch}, the result holds.\qed

\section{Tilting subcategories}

In this section, we begin with the definition of a tilting (cotilting) subcategories in an extriangulated category $\C$ with enough projectives and injectives. Then we formulate the Bazzoni  characterization for tilting (cotilting) subcategories. Finally, we generalize the Auslander-Reiten correspondence in module categories over Artin algebras in \cite{ar} to extriangulated categories.

Now we give the definition of tilting (cotilting) subcategories in an extriangulated category.
\begin{definition}\label{ti}
Let $n$ be a fixed non-negative integer, $\mathcal T$ be a subcategory of an extriangulated category $\C$ closed under direct summands. Then $\mathcal T$ is called a tilting subcategory of projective dimension $n$ if the following conditions are satisfied:
\begin{itemize}
\item[$(1)$] $\emph{pd}(\mathcal T)\leq n;$

\item[$(2)$] $\mathcal T$ is an generator for $\mathcal T^{\bot}.$

\end{itemize} We also call $\mathcal T$ a tilting subcategory for short.
$T\in\C$ is called a tilting object if $addT$ is a tilting subcategory.
\end{definition}

Dually, we define the cotilting subcategories as follows:

\begin{definition}\label{coti}
Let $n$ be a fixed non-negative integer, $\mathcal T$ be a subcategory of an extriangulated category $\C$ closed under direct summands. Then $\mathcal T$ is called a cotilting subcategory of injective dimension $n$ if the following conditions are satisfied:
\begin{itemize}
\item[$(1)$] $\emph{id}(\mathcal T)\leq n;$

\item[$(2)$] $\mathcal T$ is an cogenerator for ${}^{\bot}{\mathcal T}.$

\end{itemize} We also call $\mathcal T$ a cotilting subcategory for short.
$T\in\C$ is called a cotilting object if $addT$ is a cotilting subcategory.
\end{definition}

\begin{remark}\label{nonp}

We don't require the condition that $Proj(\C)\subset\check{\mathcal T_{n}}$ in the definition of tilting subcategories as \cite{dwzc}, we will show it is a consequence of the conditions (1), (2) in Definition \ref{ti}. Please see Remark \ref{rm}.
\end{remark}

\begin{example}\label{exa2}
\begin{itemize}
\item[$(1)$] Let $\C$ be an extriangulated category with enough projectives $Proj(\C)$ and injevtives $Inj(\C),$ then $Proj(\C)$ is a tilting subcategory and $Inj(\C)$ is a cotilting subcategory. In particular, for a triangulated category $\C$, the subcategory $0$ consisting only of zero object is a tilting subcategory and a cotilting subcategory. It is easy to see there are no other tilting (cotilting) subcategories in $\C$ except the zero category $0$. We will see for a triangulated category with a proper class of triangles in the sense of Beligiannis \cite{be}, there are some other tilting (or cotilting) subcategories (see Section 5).

\item[$(2)$] Let $R$ be an Artin algebra and $modR$ the category of finitely generated left $R-$modules and $T\in modR$ be a tilting (cotilting, resp.) module of finite projective (injective, resp.) dimension. Then $add(T)$ is a tilting (cotilting) subcategory in our sense. Actually, we will show that a $R-$module $T$ is a tilting module of finite projective dimension if and only if $T$ is a tilting object in our sense, see Remark \ref{rm}.

\item[$(3)$] Let $\mathcal E$ be an exact category with enough projectives and enough injectives. The $n$-cotilting subcategory in $\mathcal E$ defined by Haruhisa Enomoto in \cite{e} is an example of our cotilting subcategory. 
    In particular, if $R$ is an Artin algebra and $modR$ the category of finitely generated left $R-$modules, $\X\subset modR$ is a subcategory.
    Then 
    $$F_{\X}=\{0\rightarrow A\rightarrow B\rightarrow C\rightarrow0 \ exact \ in \ modR|Hom(\X,B)\rightarrow Hom(\X,C)\rightarrow0 \ is \ exact\}$$
    and
    $$F^{\X}=\{0\rightarrow A\rightarrow B\rightarrow C\rightarrow0 \ exact \ in \ modR|Hom(B,\X)\rightarrow Hom(A,\X)\rightarrow0 \ is \ exact\}$$
    defined in \cite{as} are all exact categories, hence extriangulated categories. 
    Using Remark \ref{rm} below, it is easy to see that relative tilting (cotilting, resp.) modules defined in [\cite{as}, section 3] is the same as our definition \ref{ti} and \ref{coti} in this setting.

\item[$(4)$] Let A be the path algebra of the quiver $Q:$ $$\xymatrix{3\ar[r]&2\ar[r]&1&}$$ i.e.  $A:=kQ,$ where $k$ is a field.
Then mod$A$ is an abelian category, hence an extriangulated category. Let $P(i), I(i), S(i)$ be the projective, injective and simple modules associate to the vertex $i$. It is obvious that $P(3)=I(1)$ is a porjective-injective object in modA, hence by [\cite{np}, proposition 3.30], the additive quotient categroy mod$A/P(3)$ is an extriangulated category. Since mod $A$ has almost split sequence, by [\cite{inp}, proposition 5.12], mod$A/P(3)$ also has almost split sequence and the Auslander-Reiten quiver is as follows:\\
$$\xymatrix{P(1)\ar[rd]& &S(2)\ar[rd]& &S(3) \\
 & P(2)\ar[ru]& &I(2)\ar[ru] &}$$

It is easy to verify that $P(1)\oplus P(2)$ and $P(2)\oplus S(2)$ satisfy Definition \ref{ti}, hence they are tilting objects in mod$A/P(3);$ and $I(2)\oplus S(3)$ and $I(2)\oplus S(2)$ satisfy Definition \ref{coti}, hence they are cotilting objects in mod$A/P(3).$
\end{itemize}
\end{example}

For a subcategory $\mathcal T\subset\C,$ denote by $Pres^{n}(\mathcal T)$ ($Copres^{n}(\mathcal T)$, resp.) the subcategory of all objects $A\in\C$ such that there is an $\E-$triangle sequence:
$$A'\rightarrow T_{n}\rightarrow T_{n-1}\cdots\rightarrow T_{2}\rightarrow T_{1}\rightarrow A$$
$$(A\rightarrow T_{1}\rightarrow T_{2}\cdots\rightarrow T_{n-1}\rightarrow T_{n}\rightarrow A', \mbox{resp}.)$$
in $\C$ with $T_{i}\in\mathcal T.$ It is easy to see that $\mathcal T\subset Pres^{n}(\mathcal T)$ ($\mathcal T\subset Copres^{n}(\mathcal T)$, resp.). If the equality holds, i.e. $Pres^{n}(\mathcal T)=\mathcal T$ ($Copres^{n}(\mathcal T)=\mathcal T$, resp.), we say that $\mathcal T$ is closed under $n$-image (closed under $n$-coimage, resp.).

\begin{lemma}\label{pr}
For a subcategory $\mathcal T$ of an extriangulated category $\C,$  $\emph{pd}(\mathcal T)\leq n$ iff $\mathcal T^{\bot}$ is closed under $n-$image, i.e. \ $Pres^{n}(\mathcal T^{\bot})=\mathcal T^{\bot}.$
\end{lemma}

\proof The proof for the only if part. $\mathcal T^{\bot}\subset Pres^{n}(\mathcal T^{\bot})$ is obvious from the definition. Now to show $Pres^{n}(\mathcal T^{\bot})\subset\mathcal T^{\bot}.$

For any $A\in Pres^{n}(\mathcal T^{\bot}),$ there is an $\E-$triangle sequence:
$$A'\rightarrow T_{n}\rightarrow T_{n-1}\cdots\rightarrow T_{2}\rightarrow T_{1}\rightarrow A$$
in $\C$ with $T_{i}\in\mathcal T^{\bot}.$
Hence $\E^{n+i}(T,A')=\E^{i}(T,A)$ for any $T\in\mathcal T, \ i\geq1.$ Since $\emph{pd}(\mathcal T)\leq n,$ we conclude that $\E^{i}(T,A)=0$ for any $T\in\mathcal T, \ i\geq1.$ Hence $A\in\mathcal T^{\bot},$ as desired.

The proof for the if part. Since $\C$ has enough injectives, for any $A\in\C,$ there exists an $\E-$triangle sequence:
$$A\rightarrow I_{1}\rightarrow I_{2}\cdots\rightarrow I_{n-1}\rightarrow I_{n}\rightarrow A'$$
with each $I_{i}\in Inj(\C).$ One can obtain that $A'\in\mathcal T^{\bot}$ since $Inj(\C)\subset\mathcal T^{\bot}$ and $\mathcal T^{\bot}$ is closed under $n$-image. Therefore, $\E^{n+1}(T,A)=\E(T,A')=0,$ for any $A\in\C$ and $T\in \mathcal T.$\qed

\begin{lemma}\label{pre}
Assume that  $\mathcal T$ is a self-orthonogal subcategory of an extriangulated category $\C$. Then $Pres^{n}(\mathcal T)=Pres^{n}({}_\mathcal T\X).$
\end{lemma}

\proof It is trivial that $Pres^{n}(\mathcal T)\subset Pres^{n}({}_\mathcal T\X)$ since $\mathcal T\subset {}_\mathcal T\X.$

Now to show $Pres^{n}({}_\mathcal T\X)\subset Pres^{n}(\mathcal T)$. For any $A\in Pres^{n}({}_\mathcal T\X),$ there exists an $\E-$triangle sequence:
$$A'\rightarrow X_{n}\rightarrow X_{n-1}\cdots\rightarrow X_{2}\rightarrow X_{1}\rightarrow A$$
with each $X_{i}\in {}_\mathcal T\X.$ Since $\mathcal T$ is a generator of ${}_\mathcal T\X$ and ${}_\mathcal T\X$ is closed under extensions by Lemma \ref{ortho}, then using Lemma \ref{ch}, there is an $\E-$triangle sequence:
$$V_{n}\rightarrow T_{n}\rightarrow T_{n-1}\cdots\rightarrow T_{2}\rightarrow T_{1}\rightarrow A$$
with $T_{i}\in\mathcal T.$ This shows that $A\in Pres^{n}(\mathcal T),$ hence $Pres^{n}({}_\mathcal T\X)\subset Pres^{n}(\mathcal T).$\qed

Now we prove one of our main results: the Bazzoni characterization of tilting (cotilting) categories in an extriangulated category.

\begin{theorem}\label{baz}
\begin{itemize}
\item[$(i)$]Assume that  $\mathcal T$ is a subcategory of $\C$ such that $\mathcal T$ is closed under direct summands and every object in $\mathcal T^{\bot}$ has a right $\mathcal T-$approximation. Then $\mathcal T$ is a tilting subcategory if and only if $Pres^{n}(\mathcal T)=\mathcal T^{\bot}.$
\item[$(ii)$]Assume that  $\mathcal T$ is a subcategory of $\C$ such that $\mathcal T$ is closed under direct summands and every object in ${}^{\bot}\mathcal T$ has a left $\mathcal T-$approximation. Then $\mathcal T$ is a cotilting subcategory if and only if $Copres^{n}(\mathcal T)={}^{\bot}\mathcal T.$
\end{itemize}
\end{theorem}

\proof We prove $(i),$ the proof of $(ii)$ is dually.

The only part: Using Lemma \ref{orhton}, one has $\mathcal T^{\bot}={}_\mathcal T\X\subset Pres^{n}(\mathcal T).$ We only need to show $Pres^{n}(\mathcal T)\subset\mathcal T^{\bot}.$ For any $A\in Pres^{n}(\mathcal T),$ there is an $\E-$triangle sequence:
$$A'\rightarrow T_{n}\rightarrow T_{n-1}\cdots\rightarrow T_{2}\rightarrow T_{1}\rightarrow A$$
with each $T_{i}\in \mathcal T.$
Hence $\E^{n+i}(T,A')=\E^{i}(T,A)$ for any $T\in\mathcal T, \ i\geq1.$ Since $\emph{pd}(\mathcal T)\leq n,$ we have that $\E^{i}(T,A)=0$ for any $T\in\mathcal T, \ i\geq1.$ Therefore $A\in\mathcal T^{\bot}$ and $Pres^{n}(\mathcal T)\subset\mathcal T^{\bot}.$

The if part: Since $\mathcal T\subset Pres^{n}(\mathcal T)=\mathcal T^{\bot},$ one can see that $\mathcal T$ is a self-orthonogal subcategory of $\C.$

We first show that $\mathcal T$ is a generator of $\mathcal T^{\bot}.$
For any object $A\in\mathcal T^{\bot}=Pres^{n}(\mathcal T),$ there is an $\E-$triangle $\xymatrix{A'\ar[r]^{x}&T_{0}\ar[r]^{y}&A\ar@{-->}[r]^{\delta}&}$ with $T_{0}\in\mathcal T.$
Since $\mathcal T^{\bot}$ has a right $\mathcal T-$approximation, one has a right $\mathcal T-$approximation $f\colon T_{1}\to A$ of $A$ with $T_{1}\in\mathcal T.$
Then we have the following commutative diagram:
\begin{equation}\label{t2}
\begin{array}{l}
\xymatrix{A^{'}\ar[r]^{a}\ar@{=}[d]&K\ar[r]^{b}\ar[d]^{g}&T_{1} \ar@{-->}[r]^{f^{*}\del}\ar[d]^{f}&\\
A^{'}\ar[r]^{x}&T_{0}\ar[r]^{y}&A\ar@{-->}[r]^{\del}&.}
\end{array}
\end{equation}
According to the dual result of Lemma \ref{pull}, we have the following $\E$-triangle
$$\xymatrix{K\ar[r]^{\binom{g}{b}\quad}&T_{0}\oplus T_{1}\ar[r]^{\quad(y,\ -f)}&A\ar@{-->}[r]^{a_{*}\delta}&}$$
with$\quad(y,\ -f)$ a right $\mathcal T-$approximation of $A.$
By Lemma \ref{lon} and Lemma \ref{long}, for any $T\in\mathcal T,$ one obtains the following long exact sequence:
$$\C(T, T_{0}\oplus T_{1})\xrightarrow{\C(T,(y,\, -f))}\C(T, A)\xrightarrow{a_{*}\delta^{\sharp}_T}
\E(T, K)\xrightarrow{\E(T,\binom{g}{b})}\E(T, T_{0}\oplus T_{1})\rightarrow$$ $$\cdots\rightarrow\E^{i}(T, K)\xrightarrow{\E^{i}(T,\binom{g}{b})}\E^{i}(T, T_{0}\oplus T_{1})
\xrightarrow{\E^{i}(T,(y,\, -f))}\E^{i}(T, A)\rightarrow\cdots$$
Since $\C(T,(y,\ -f))$ is an epimorphism and $\E^{i}(T, T_{0}\oplus T_{1})=\E^{i}(T, A)=0, \forall i\geq1,$ one can conclude that $\E^{i}(T,K)=0$ for any $T\in\mathcal T$ and any $i\geq1,$ i.e. $K\in \mathcal T^{\bot}.$ Hence $\mathcal T$ is a generator of $\mathcal T^{\bot}.$

Now to show $\emph{pd}(\mathcal T)\leq n.$
Applying Lemma \ref{orthon}, one has $\mathcal T^{\bot}={}_\mathcal T\X.$
Hence, the following equality holds
$$Pres^{n}(\mathcal T^{\bot})=Pres^{n}({}_\mathcal T\X)=Pres^{n}(\mathcal T)=\mathcal T^{\bot}$$
by Lemma \ref{pre}. Using Lemma \ref{pr}, one obtains that $\emph{pd}(\mathcal T)\leq n.$\qed

\begin{remark}\label{bazz}
One can see that in the proof of this theorem, we do not use the fact that $\C$ is a Krull-Schmit category. Hence the theorem holds for any extriangulated category with enough projectives and injectives. A typical example is the characterization of $n$-tilting and $n$-cotilting modules in Mod$R$ \cite{ba} where $R$ is a ring and Mod$R$ is the category of left $R-$modules.
\end{remark}

In the rest of this section, we will establish the Auslander-Reiten correspondence between tilting (cotilting, resp.) subcategories and coresolving covariantly finite (resloving contravariantly finite, resp.) subcategories satisfying some cogenerating (generating, resp.) conditions.

A subcategory $\X\subset\C$ is called coresolving if it contains $Inj(\C),$ closed under extensions and Cone of deflations. Resolving subcategory can be defined dually.

\begin{lemma}\label{core}
Assume that  $\X$ is a coresolving subcatefory of $\C,$  denote $\sideset{^{\bot_{1}}}{}{\mathop{\X}}$ by $\Y$. Then $\Y=\sideset{^{\bot}}{}{\mathop{\X}}.$
\end{lemma}

\proof It is easy to see that $\sideset{^{\bot}}{}{\mathop{\X}}\subset\Y.$ On the other hand, For any $X\in\X,$ there is an injective resolution of $X$:
$$X\xrightarrow{d_{0}}I_{0}\xrightarrow{d_{1}}I_{1}\rightarrow\cdots  \rightarrow I_{n-1}\xrightarrow{d_{n}}I_{n}\rightarrow\cdots $$
with each $I_{i}\in Inj(\C).$ Since $Inj(\C)\subset\X$ and $\X$ is a coresolving subcategory, we have $\Sigma^{i}X\in\X,\ \forall i\geq0.$ Hence for any $Y\in\Y,$ $\E^{i}(Y,X)=\E(Y,\Sigma^{i-1}X)=0,\ \forall i\geq1.$\qed

\begin{lemma}\label{cores}
Assume that  $\X$ is an coresolving covariantly finite subcategory of $\C.$ Let $\mathcal T=\sideset{^{\bot}}{}{\mathop{\X}}\cap\X.$ Then $\mathcal T$ is a generator of $\X.$
\end{lemma}

\proof According to Lemma \ref{core}, one has $\sideset{^{\bot_{1}}}{}{\mathop{\X}}=\sideset{^{\bot}}{}{\mathop{\X}}.$
Using Corollary \ref{cocon}, for any $X\in\X,$ there is an $\E-$triangle:
$$\xymatrix{X'\ar[r]^{x}&T\ar[r]^{y}&X\ar@{-->}[r]^{\delta}&}$$
with $X'\in\X$ and $T\in\sideset{^{\bot}}{}{\mathop{\X}}.$
As $\X$ is closed under extensions, then $T\in\X.$
One has $T\in\mathcal T,$ hence $\mathcal T$ is a generator of $\X.$\qed

We say that $\mathcal T\subset\X$ is an Extprojective generator of $\X$ if $\mathcal T$ is a generator of $\X$ and $\E(T,X)=0, \forall \ T\in\mathcal T, X\in\X.$ Now we state the following result, which is an extriangulated category version of the dual result of [\cite{ab}, Proposition 3.6]. We omit the proof:

\begin{lemma}\label{abc}
Assume that $\mathcal T\subset\X$ are subcategories of $\C$ such that $\mathcal T$ and $\X$ are closed under direct summands and $\mathcal T$ is a Extprojective generator of $\X.$ Then we have $\check{\mathcal T}_{n}=\check{\X_{n}}\cap\sideset{^{\bot}}{}{\mathop{\X}}.$
\end{lemma}

Now we introduce a partial order on tilting subcategories as a generalization of the partial order on tilting modules introduced by Happel and Unger \cite{hu}.

\begin{definition}\label{ord}
Assume that $\mathcal T_{i}, \ i=1, 2$ are tilting subcategories of $\C,$ we write $\mathcal T_{1}\leq\mathcal T_{2}$ if $\mathcal T_{1}\subset\mathcal T_{2}^{\bot}.$
\end{definition}

The following lemma ensure that the relation defined above on tilting subcategories is a partial order.

\begin{lemma}\label{orde}
Assume that $\mathcal T, \mathcal T_{i}, \ i=1, 2$ are tilting subcategories of $\C.$ We have:
\begin{itemize}
\item[$(i)$] $\mathcal T=\mathcal{T}^{\bot}\cap {}^{\bot}{(\mathcal{T}^{\bot})};$
\item[$(ii)$] $\mathcal T_{1}\leq\mathcal T_{2}$ if and only if $\mathcal T_{1}^{\bot}\subset\mathcal T_{2}^{\bot}$
\end{itemize}
\end{lemma}

\proof $(i)$ It is obvious that $\mathcal T\subset\mathcal{T}^{\bot}\cap {}^{\bot}{(\mathcal{T}^{\bot})}$ since $\mathcal T$ is self-orthonogal. On the other hand, since $\mathcal T$ is a generator of $\mathcal{T}^{\bot},$ for any $A\in\mathcal{T}^{\bot}\cap {}^{\bot}{(\mathcal{T}^{\bot})},$ there is an $\E-$triangle  \xymatrix{A'\ar[r]^{x}&T\ar[r]^{y}&A\ar@{-->}[r]^{\delta}&} with $A'\in \mathcal{T}^{\bot}$ and $T\in\mathcal T.$ Hence the $\E-$triangle splits, we have $A\in\mathcal T$ as $\mathcal T$ is closed under direct summands.

$(ii)$ The if part:  It is easy to see that $\mathcal T_{1}\leq\mathcal T_{2}$ since $\mathcal T_{1}\subset\mathcal T_{1}^{\bot}.$

The only if part: Let $A\in\mathcal T_{1}^{\bot},\ T\in\mathcal T_{2}.$ Since $\mathcal T_{1}$ is a generator of $\mathcal T_{1}^{\bot},$ we have the following $\E-$triangle sequence:
$$\cdots\rightarrow T_{n}\xrightarrow{d_{n}}T_{n-1}\rightarrow\cdots  \rightarrow T_{1}\xrightarrow{d_{1}}T_{0}\xrightarrow{d_{0}}A$$
in $\C$ with $T_{i}\in\mathcal T_{1}\subset\mathcal T_{2}^{\bot}, d_{n}=g_{n-1}f_{n},$ where $\xymatrix{A_{n}\ar[r]^{g_{n-1}}&T_{n-1}\ar[r]^{f_{n-1}}&A_{n-1}\ar@{-->}[r]^{\delta_{n-1}}&}$ is an $\E-$triangle and $A_{n}\in\mathcal T_{1}.$ Hence, we have $\E^{n+i}(T,A_{n})=\E^{i}(T,A)=0,$ for any $i\geq1$ since $\emph{pd}(\mathcal T_{2})\leq n.$ Then we have $\mathcal T_{1}^{\bot}\subset\mathcal T_{2}^{\bot},$ this finishs the proof.\qed

\begin{corollary}\label{order}
The relation defined on tilting subcategories is a partial order.
\end{corollary}

\proof $\mathcal T\leq\mathcal T$ since $\mathcal T$ is self-orthonogal.

If $\mathcal T_{1}\leq\mathcal T_{2}$ and $\mathcal T_{2}\leq\mathcal T_{3},$ then $\mathcal T_{1}^{\bot}\subset\mathcal T_{2}^{\bot}\subset\mathcal T_{3}^{\bot}$ by Lemma \ref{orde}, hence $\mathcal T_{1}\leq\mathcal T_{3}.$

If $\mathcal T_{1}\leq\mathcal T_{2}$ and $\mathcal T_{2}\leq\mathcal T_{1},$ then $\mathcal T_{1}^{\bot}=\mathcal T_{2}^{\bot}$ by Lemma \ref{orde}, hence $\mathcal T_{1}=\mathcal T_{2}$\qed

Now we prove the Auslander-Reiten correspondence for tilting subcategories in an extriangulated category.
\begin{theorem}\label{arc} Let $\C$ be an extriangulated category with enough projectives and enough injectives. Then
\begin{itemize}
\item[$(i)$]The assignments $\mathcal T\mapsto \mathcal T^{\bot}$ and $\X\mapsto \sideset{^{\bot}}{}{\mathop{\X}}\cap\X$ give a one-to-one correspondence between the class of tilting subcategories and coresolving covariantly finite subcategories $\X,$ which are closed under summands with $\check{\X_{n}}=\C.$
\item[$(ii)$]The assignments $\mathcal T\mapsto {}^{\bot}\mathcal T$ and $\X\mapsto \X^{\bot}\cap\X$ give a one-to-one correspondence between the class of cotilting subcategories and resolving contravariantly finite subcategories $\X,$ which are closed under summands with $\hat{\X_{n}}=\C.$
\end{itemize}
\end{theorem}

\proof We only prove $(i),$ the proof of $(ii)$ is dually.

(1) Assume that $\mathcal T$ is a tilting subcategory. Let $\X=\mathcal T^{\bot}.$
For any element $A\in\C,$ there is an $\E-$triangle sequence:
$$A\xrightarrow{d_{0}}I_{0}\xrightarrow{d_{1}}I_{1}\rightarrow\cdots  \rightarrow I_{n-2}\xrightarrow{d_{n-1}}I_{n-1}\xrightarrow{d_{n-1}}A'$$
with $I_{i}\in Inj(\C)\subset\X.$
One has $\E^{n+i}(T,A)=\E^{i}(T,A')=0, \forall i\geq1$ and for any $T\in\mathcal T$ since $\emph{pd}(\mathcal T)\leq n.$ This implies that $A'\in\mathcal T^{\bot}=\X,$ hence $A\in\check{\X_{n}}.$

It is obvious that $\X$ is coresolving, we only need to show that $\X$ is covariantly finite.
Since $\mathcal T$ is a generator of $\X$ and $\X$ is closed under extensions with $\check{\X_{n}}=\C,$ using Lemma \ref{chang}, one can see for any object $A\in\C,$ there is an $\E-$triangle $\xymatrix{U\ar[r]^{x}&V\ar[r]^{y}&A\ar@{-->}[r]^{\delta}&}$ in $\C$ with $U\in\X$ and $V\in\check{\mathcal T_{n}}.$ Using {\rm (ET4)}, we have the following commutative diagram:
$$\xymatrix{U\ar[r]^{x}\ar@{=}[d]&V\ar[r]^{y}\ar[d]^{a}&A\ar[d]^{f}\\
U\ar[r]^{m}&T\ar[d]^{b}\ar[r]^{n}&A_{1}\ar[d]^{g}\\
&A_{2}\ar@{=}[r]&A_{2}}$$
with $T\in\mathcal T, A_{2}\in\check{\mathcal T}_{n-1}.$
One obtain that $A_{1}\in\X$ since $U$ and $T$ are in $\X$ and $\X$ is coresolving.
We state that $f \colon A\to A_{1}$ is a left $\X-$approximation of $A$ by Lemma \ref{ort}, hence $\X$ is covariantly finite. Using Lemma \ref{orde}, we also have that $\mathcal T=\mathcal{T}^{\bot}\cap {}^{\bot}{(\mathcal{T}^{\bot})}=\X\cap {}^{\bot}{\X}.$

(2) Suppose that $\X$ is a coresolving covariantly finite subcategory, which is closed under summands with $\check{\X_{n}}=\C.$ Let $\mathcal T=\sideset{^{\bot}}{}{\mathop{\X}}\cap\X.$ Now to show that $\emph{pd}(\mathcal T)\leq n.$ For any object $A\in\C=\check{\X_{n}},$ there is an $\E-$triangle sequence:
$$A\rightarrow X_{0}\rightarrow X_{1}\cdots\rightarrow X_{n-1}\rightarrow X_{n}$$
with $X_{i}\in\X.$ It follows that $\E^{n+1}(T,A)=\E(T,X_{n})=0$ for any $T\in\mathcal T,$ hence $\emph{pd}(\mathcal T)\leq n$ by Lemma \ref{projec}.

By Lemma \ref{cores}, $\mathcal T$ is a generator of $\X.$ Using the same argument as in (1), it is easy to see that for any object $A\in\mathcal T^{\bot},$ there is an $\E-$triangle $\xymatrix{A\ar[r]^{f}&A_{1}\ar[r]^{g}&A_{2}\ar@{-->}[r]^{\delta}&}$ with $A_{1}\in\X$ and $A_{2}\in\check{\mathcal T}_{n-1},$ which splits by Lemma \ref{ort}, hence $A\in\X.$ So one can see that $\X=\mathcal T^{\bot},$ i.e. $\X=(\sideset{^{\bot}}{}{\mathop{\X}}\cap\X)^{\bot}=\mathcal{T}^{\bot}.$ It follows that $\mathcal T$ is a generator of $\mathcal T^{\bot}.$
This completes the proof.\qed

\begin{remark}\label{rm}
\item[$(1)$] It is easy to see that if $\mathcal T$ is a tilting subcategory of $\C$, then $\mathcal T$ is an Extprojective generator of $\X=\mathcal T^{\bot}.$ Using Lemma \ref{abc} and Theorem \ref{arc}, we have $Proj(\C)\subset\check{\mathcal T_{n}}$ since $Proj(\C)\subset\sideset{^{\bot}}{}{\mathop{\X}}$ and $\check{\X_{n}}=\C.$ Similarly, one can see that  $Inj(\C)\subset\hat{\mathcal T_{n}}$ if $\mathcal T$ is a cotilting subcategory of $\C$. Hence if $\C=modR,$ where $R$ is an Artin algebra and $modR$ is the category of finitely generated left $R-$modules, then our definition \ref{ti} and \ref{coti} coincide with the usual definition of (co)tilting objects of finite projective (injective) dimension in $modR.$
\item[$(2)$]Let $R$ be a Artin algebra and mod$R$ denote the category of finitely generated left $R-$modules. Then mod$R$ is an abelian category with enough projectives and injectives. Applying Theorem \ref{arc} to mod$R$, one can deduce the classical Auslander-Reiten correspondence {[\cite{ar}, Theorem 5.5]}.
\end{remark}

\section{Applications}

In this section, we give some applications of our results. We first show that the relative tilting subcategories in an abelian category $\A$ with a complete cotorsion triple $(\X, \mathcal Z, \Y)$ (for more detail, see \cite{dwzc}) are tilting subcategories in our sense. Therefore the results in \cite{dwzc} follow from our main results above. We also show that for a triangulated category $\C$ with a proper class of triangles $\mathcal{E}$ introduced by Beligiannis \cite{be} form an extriangulated category. If in addition $\mathcal{E}$ has enough projectives and injectives, we obtain the corresponding results about tilting (cotilting) subcategories in this setting.

\subsection{Abelian categories with a complete cotorsion triple}

In this subsection, $\A$ is a Krull-Schmit abelian category with enough projectives $\mathcal{P}$ and injectives $\mathcal{I}$. All subcategories are full additive subcategories of $\A$ closed under isomorphisms. A pair $(\X, \Y)$ of subcategories in $\A$ is called a cotorsion pair if $\X=\sideset{^{\bot1}}{}{\mathop{\Y}}$ and $\Y=\X^{\bot1}.$ The cotorsion pair $(\X, \Y)$ is said to be complete provided that, for any $A\in \A,$ there are short exact sequences $0\rightarrow Y\rightarrow X\rightarrow A\rightarrow0$ and $0\rightarrow A\rightarrow Y'\rightarrow X'\rightarrow0$ in $\A$ with $X, X'\in \X$ and $Y, Y'\in\Y.$ The cotorsion pair $(\X, \Y)$ is said to be hereditary if $\X$ is resolving, or equivalently if $\Y$ is coresolving.

Let $\X, \Y, \mathcal Z$ be subcategories of $\A.$ The triple $(\X, \mathcal Z, \Y)$ is said to be a cotorsion triple (see \cite{c}) if both $(\X, \mathcal Z)$ and $(\mathcal Z, \Y)$ are cotorsion pairs.
A cotorsion triple is said to be complete ( hereditary, resp.) if both $(\X, \mathcal Z)$ and $(\mathcal Z, \Y)$ are complete (hereditary, resp.).
A standard example of cotorsion triple is $(\mathcal P, \mbox{mod}R, \mathcal I)$ where $R$ is a Artin algabra and $\mathcal P, \mathcal I$ are subcategories of projecitves and injectives of mod$R,$ resp.

Given a complete hereditary cotorsion triple $(\X, \mathcal Z, \Y)$ of $\A,$ one can define an exact structure on $\A$ as follows:
Let $F_{\X}$ be the set of exact sequences of $\A$ which are also $Hom_{\A}(\X,-)$ exact, i.e.
$$F_{\X}=\{0\rightarrow A\rightarrow B\rightarrow C\rightarrow0 \ exact \ in \ \A|Hom(\X,B)\rightarrow Hom(\X,C)\rightarrow0 \ is \ exact\}.$$
Then by [\cite{drssk}, Proposition 1.7], $(\A, F_{\X})$ is an exact category, hence an extriangulated category. We denote the exact category $(\A, F_{\X})$ by $\mathcal{E},$ and the subcategories of  projective objects or of injective objects in $\mathcal{E}$ by $\mathcal{P(\mathcal{E})}$ or $\mathcal{I(\mathcal{E})}$, respectively.

A complex $\bm{C}$ in $\A$ is called $\X-$ exact if it is exact in $\A$ and $Hom_{\A}(\X,-)$ exact.
Since $(\X, \mathcal Z)$ is a complete cotorsion pair, for any object $A\in\A,$ there is an proper $\X-$resolution $\bm{X}$ of $A$:
$$\dots\rightarrow X_{1}\rightarrow X_{0}\rightarrow A\rightarrow0$$
such that $X_{i}\in\X$ and $\bm{X}$ is $Hom_{\A}(\X,-)$ exact.

The $\X-$projective dimension of $A$ is the least non-negative number $n$ such that there is an exact sequence $0\rightarrow X_{n}\dots\rightarrow X_{1}\rightarrow X_{0}\rightarrow A\rightarrow0$ in $\A$ with $X_{i}\in\X.$ We denote the $\X-$projective dimension of $A$ by $\X-\emph{pd}(A).$

The proper $\X-$resolution of $A\in\A$ is unique up to homotopy.
Hence, for any $B\in\A,$ one can define the relative $\X$ cohomology group
$$Ext_{\X}^{i}(M,N)=H_{-i}(Hom_{\A}(\bm{X},N)).$$
Similarly, the terms of proper $\Y-$resolution $\bm{Y},$ $\Y-$injective dimension and relative $\Y$ cohomology group $Ext_{\Y}^{i}(-,-)$ are defined dually.

\begin{lemma}{[\cite{dwzc}, Proposition 2.6]}\label{sym}
Let $(\X, \mathcal Z, \Y)$ be a complete hereditary cotorsion triple in $\A.$ Then, for all objects $A, B\in\A$ and all $i\in\mathbb{Z},$ there exist isomorphisms
$$Ext_{\X}^{i}(A,B)=Ext_{\Y}^{i}(A,B).$$
\end{lemma}

According to the arguments above, we have the following result.
\begin{lemma}\label{pi}
$\mathcal{E}$ is an exact category with enough projectives and injectives. In particular, $\mathcal{P(\mathcal{E})}=\X$ and $\mathcal{I(\mathcal{E})}=\Y$.
\end{lemma}

\proof Firstly, to show that $\mathcal{E}$ has enough projectives with $\mathcal{P(\mathcal{E})}=\X.$ For any $\X-$exact sequence $0\rightarrow A\rightarrow B\rightarrow X\rightarrow0$ with $X\in\X,$ the sequence $0\rightarrow Hom(X,A)\rightarrow Hom(X,B)\rightarrow Hom(X,X)\rightarrow0$ is exact. Hence the $\X-$exact sequence splits, i.e. $Ext_{\X}(X,A)=0,$ for any $X\in\X, \ A\in\A.$ Then we have $\X\subset\mathcal{P(\mathcal{E})}.$

Since $(\X, \mathcal Z)$ is a complete cotorsion pair, for any $A\in\A,$ there is an $\X-$exact sequence $0\rightarrow Z\rightarrow X\rightarrow A\rightarrow0$ with $X\in\X, \ Z\in\mathcal{Z}.$ Hence $\mathcal{E}$ has enough projectives. In particular, if $A\in\mathcal{P(\mathcal{E})},$ it is obvious that the sequence splits, one has $A\in\X,$ hence $\mathcal{P(\mathcal{E})}=\X.$

Secondly, to show that $\mathcal{E}$ has enough injectives with $\mathcal{I(\mathcal{E})}=\Y.$ Using Lemma 5.1, one can see that $Ext_{\X}(A,Y)=Ext_{\Y}(A,Y)=0$ for any $Y\in\Y,$ hence, $\Y\subset\mathcal{I(\mathcal{E})}.$

Since $(\mathcal Z, \Y)$ is a complete cotorsion pair,  for any $A\in\A,$ there is a $\Y-$exact sequence $ 0\rightarrow A\xrightarrow{f} Y\xrightarrow{g} Z\rightarrow0$ with $Y\in\Y, \ Z\in\mathcal{Z}.$ Using the same method as that in the proof of Lemma 2.5 in \cite{dwzc}, one can see that the sequence is also $\X-$exact. Therefore, $\mathcal{E}$ has enough injectives. In particular, if $A\in\mathcal{I(\mathcal{E})},$ the sequence splits, hence $A\in\Y.$ We have that $\mathcal{I(\mathcal{E})}=\Y.$ This completes the proof.\qed

Given a complete hereditary cotorsion triple $(\X, \mathcal Z, \Y)$ in an abelian category $\mathcal A$ with enough projectives and injectives, we have shown that the associate category $\mathcal{E}$ is an exact category with enough projectives and injectives.

A subcategory $\mathcal M$ is said to be $\X-$self-orhtonogal if $Ext_{\X}^{i\geq1}(M,M')=0$ for any $M, M'\in \mathcal M.$ For an $\X-$self-orhtonogal subcatgory $\mathcal M$, we can define ${}_{\X}{\hat{\mathcal M}}_{n},$ ${}_\X{\check{\mathcal M}}_{n},$ ${}_{\X}{\hat{\mathcal M}},$ ${}_\X{\check{\mathcal M}},$ $\mathcal M^{\X\bot}$ and ${}^{\X\bot}{\mathcal M}$ as in [DWZC, section 3].

Recall the definition of $n$-$\X-$tilting subcategories defined in [DWZC, Definition 4.1] as follows:

\begin{definition}\label{rti}{\emph{[\cite{dwzc}, Definition 4.1]}}
Let $n$ be a positive integer and $\mathcal M$ be a subcategory of $\A$ closed under direct summands. Then $\mathcal M$ is called a  $n-\X-$tilting (with respect to $(\X, \mathcal Z, \Y)$) if the following conditions are satisfied:
\begin{itemize}
\item[$(1)$] $\X-\emph{pd}(\mathcal M)\leq n;$

\item[$(2)$] $\mathcal M$ is an $\X$generator for $\mathcal M^{\X\bot};$

\item[$(3)$] $\X\subset {}_\X{\check{\mathcal M}_{n}}.$
\end{itemize}
\end{definition}

According to Remark \ref{rm} and Lemma \ref{pi}, one can see that the above definition is the same as Definition \ref{ti} when consider $\mathcal{E}$ as an extriangulated category with enough projectives and injectives.

Using Theorem \ref{baz} and Theorem \ref{arc} we have the following Bazzoni type characterization of $n$-$\X-$tilting subcategories and the Auslander-Reiten correspondence in \cite{dwzc}.

\begin{corollary}\label{dbaz}{\emph{[\cite{dwzc}, Theorem 5.5]}}
Assume that  $\mathcal M\subset\C$ such that $\mathcal M$ is closed under direct summands and every object in $\mathcal M^{\X\bot}$ has a right $\mathcal M-$approximation. Then $\mathcal M$ is a n-$\X-$tilting subcategory (with respect to $(\X, \mathcal Z, \Y)$) iff $Pres^{n}(\mathcal T)=\mathcal T^{\bot}.$
\end{corollary}

\begin{corollary}\label{darc}{\emph{[\cite{dwzc}, Theorem 7.7]}}
The assignments $\mathcal M\mapsto \mathcal T^{\bot}$ and $\X\mapsto \sideset{^{\bot}}{}{\mathop{\mathcal{H}}}\cap\mathcal{H}$ induce a one-to-one correspondence between the class of $n-\X-tilting$ subcategories (with respect to $(\X, \mathcal Z, \Y)$) and subcategories $\mathcal{H}$ of $\A,$ which are $\X-$resolving covariantly finite and closed under direct summands, such that $\check{\mathcal{H}}_{n}=\A.$
\end{corollary}

\subsection{Triangulated categories with a proper classes of triangles.}
Beligiannis considered in \cite{be} the homological theory of a triangulated category $\C$ by using the proper class $\mathcal{E}$ of triangles.
We will show that the triangulated category $\C$ together with the proper class $\mathcal{E}$ forms an extriangulated category, denoted by $(\C, \mathcal{E})$.
If in addition, it has enough projectives and injectives, we can define tilting subcategories in $(\C, \mathcal{E})$ as Definition \ref{ti} and then get the Bazonni characterization of tilting subcategories and the Auslander-Reiten correspondence in this situation.

Let $(\C,\Sigma,\vartriangle)$ be a Krull-Schmit triangulated category, $\Sigma$ an automorphism of $\C$ called suspension and $\vartriangle$ a class of diagrams in $\C$ of the form $A\rightarrow B\rightarrow C\rightarrow\Sigma A,$ called distinguished triangles, satisfying axioms (Tr1)-(Tr4)[V].

A triangle $A\xrightarrow{f} B\xrightarrow{g} C\xrightarrow{h} \Sigma A$ is called splitting if $h=0.$ We denote the full subcategory of $\Sigma$ consisting of splitting triangles by $\vartriangle_{0}.$

\begin{definition}{\cite{be}}\label{pt}
Let $\C$ be a triangulated category. Assume that $\mathcal{E}\subset\vartriangle$ is a class of triangles.
\begin{itemize}
\item[{\rm (i)}] $\mathcal{E}$ is closed under base change if for any triangle $A\xrightarrow{f} B\xrightarrow{g} C\xrightarrow{h} \Sigma A\in\mathcal{E}$ and any morphism $\varepsilon\colon E\to C$ as in the diagram of [\cite{be}, Proposition 2.1(i)], the triangle $A\xrightarrow{f'} G\xrightarrow{g'} E\xrightarrow{h\varepsilon} \Sigma A\in\mathcal{E}.$
\item[{\rm (ii)}] $\mathcal{E}$ is closed under cobase change if for any triangle $A\xrightarrow{f} B\xrightarrow{g} C\xrightarrow{h} \Sigma A\in\mathcal{E}$ and any morphism $\alpha\colon A\to D$ as in the diagram of [\cite{be}, Proposition 2.1(ii)], the triangle $D\xrightarrow{f'} F\xrightarrow{g'} C\xrightarrow{\Sigma \alpha h} \Sigma D\in\mathcal{E}.$
\item[{\rm (iii)}] $\mathcal{E}$ is closed under suspension if for any triangle $A\xrightarrow{f} B\xrightarrow{g} C\xrightarrow{h} \Sigma A\in\mathcal{E}$ and any $i\in\mathbb{Z},$ the triangle $\Sigma A\xrightarrow{(-1)^{i}\Sigma^{i}f} B\xrightarrow{(-1)^{i}\Sigma^{i}g} C\xrightarrow{(-1)^{i}\Sigma^{i}h} \Sigma^{i+1} A\in\mathcal{E}.$
\item[{\rm (i)}] $\mathcal{E}$ is saturated if in the diagram [\cite{be}, Proposition 2.1(i)], the third vertical and the second horizontal triangles are in $\mathcal{E}$ imply  that the triangle $A\xrightarrow{f} B\xrightarrow{g} C\xrightarrow{h} \Sigma A\in\mathcal{E}.$
\end{itemize}
\end{definition}

\begin{definition}\label{prt}{\cite{be}}
A full subcategory $\mathcal{E}\subset\vartriangle$ is called a proper class of triangles if the following conditions hold.
\begin{itemize}
\item[{\rm (i)}] $\mathcal{E}$ is closed under isomorphisms, finite coproducts and $\vartriangle_{0}\subset\mathcal{E}\subset\vartriangle.$
\item[{\rm (ii)}] $\mathcal{E}$ is closed under suspensions and is saturated.
\item[{\rm (iii)}] $\mathcal{E}$ is closed under base change and cobase change.
\end{itemize}
\end{definition}

We can define the $\mathcal{E}-$phantom map for a proper class of triangles in $\C,$ which can be used to construct an additive bifunctor.

\begin{definition}\label{ph}{\cite{be}}
Let $\mathcal{E}$ be a proper class of triangles in $\C$ and let $A\xrightarrow{f} B\xrightarrow{g} C\xrightarrow{h} \Sigma A$ be a triangle in $\mathcal{E}.$
\begin{itemize}
\item[{\rm (i)}] The morphism $f\colon A\to B$ is called an $\mathcal{E}-$proper monic.
\item[{\rm (ii)}] The morphism $g\colon B\to C$ is called an $\mathcal{E}-$proper epic.
\item[{\rm (iii)}] The morphism $h\colon C\to \Sigma A$ is called an $\mathcal{E}-$phantom map.
\end{itemize}
Let $Ph_{\mathcal{E}}(A,B)$ be the set of all $\mathcal{E}-$phantom maps from $A$ to $B$, the class of $\mathcal{E}-$phantom maps is denoted by $Ph_{\mathcal{E}}(\C).$
\end{definition}

Let $\mathcal{E}$ be a proper class of triangles in $\C.$ Given $A, C\in \C$ and consider the class $\mathcal{E}^{*}(C,A)$ of all triangles $A\xrightarrow{f} B\xrightarrow{g} C\xrightarrow{h} \Sigma A$ in $\mathcal{E}.$ One can define a relation in $\mathcal{E}^{*}(C,A)$ as follows. If ($T_{i}$): $A\xrightarrow{f_{i}} B_{i}\xrightarrow{g_{i}} C\xrightarrow{h_{i}} \Sigma A,$ $i=1,2,$ are elements of $\mathcal{E}^{*}(C,A)$, then define $(T_{1})\sim (T_{2})$ if there is a morphism of triangles
$$\xymatrix{A\ar[r]^{f_{1}}\ar@{=}[d]&B_{1}\ar[r]^{g_{1}}\ar[d]^{\alpha}&C\ar[r]^{h_{1}}\ar@{=}[d]& \Sigma A\ar@{=}[d]\\
A\ar[r]^{f_{2}}&B_{2}\ar[r]^{g_{2}}&C\ar[r]^{h_{2}}&
\Sigma A}$$

It is obvious that $\alpha$ is an isomorphism and $\sim$ is an equivalence relation on the class $\mathcal{E}^{*}(C,A).$ Using base change and cobase change, one can define a sum in the class $\mathcal{E}(C,A):= \mathcal{E}^{*}(C,A)/\sim.$ It is easy to see that $\mathcal{E}(-,-): \C^{op}\times\C\rightarrow\A b$ is an additive bifunctor.
According to [Be, Corollary 2.6], one has $\mathcal{E}(-,-)\cong Ph_{\mathcal{E}}(-,\Sigma-),$ hence $Ph_{\mathcal{E}}(-,\Sigma-)$ defines a bifunctor on $\C.$
For any $A, C\in\C$ and $h\in Ph_{\mathcal{E}}(C,\Sigma A),$ one can consider the triangle $A\xrightarrow{f} B\xrightarrow{g} C\xrightarrow{h} \Sigma  A$ in $\mathcal{E}$ as a realization of $h$ since $\mathcal{E}$ is closed under base change and cobase change, and we denote it by $\s(h)$. In this way $Ph_{\mathcal{E}}(-,\Sigma-)$ can be considered as the additive bifunctor $\E$ and every equivalence class of triangles in $\mathcal{E}$ can be considered as an $\E-$triangle. Under this circumstance, we have:

\begin{lemma}\label{pe}
$(\C,Ph_{\mathcal{E}},\s)$ is an extriangulated category.
\end{lemma}

\proof It is easy to verify that (ET1)-(ET3) are satisfied for the class of $\E-$triangles.

Now to show (ET4). Let $A\xrightarrow{f} B\xrightarrow{g} C\rightarrow \Sigma A\in\mathcal{E}$ and $B\xrightarrow{a} C\xrightarrow{b} F\rightarrow \Sigma A\in\mathcal{E}$ are $\E-$triangles in $(\C,\mathcal{E}).$ Using (Tr4), we have the following commutative diagram:
$$\xymatrix{A\ar[r]^{f}\ar@{=}[d]&B\ar[r]^{g}\ar[d]^{a}&C\ar[d]^{c}\\
A\ar[r]^{af}&C\ar[d]^{b}\ar[r]^{m}&E\ar[d]^{d}\\
&F\ar@{=}[r]&F}$$
Since $B\xrightarrow{a} C\xrightarrow{b} F\rightarrow \Sigma A\in\mathcal{E}$ is an $\E-$triangle in $(\C,\mathcal{E})$ and $\mathcal{E}$ is closed under cobase change, we have that $C\xrightarrow{c} E\xrightarrow{d} F\rightarrow \Sigma C\in\mathcal{E}.$ Therefore, $A\xrightarrow{af} C\xrightarrow{m} E\rightarrow \Sigma A\in\mathcal{E}$ as $A\xrightarrow{f} B\xrightarrow{g} C\xrightarrow{h} \Sigma A\in\mathcal{E}$ and $\mathcal{E}$ is saturated. Hence, (ET4)holds. Similarly, one can show that (ET4$\op$) holds. Hence, $(\C,Ph_{\mathcal{E}},\s)$ is an extriangulated category.\qed

If in addition, $(\C,Ph_{\mathcal{E}},\s)$ has enough projectives and injectives, we can define tilting subcategories in $(\C,Ph_{\mathcal{E}},\s)$ as follows:

\begin{definition}\label{ptti}
Let $\mathcal T$ be a subcategory of $\C$ closed under direct summands. Then $\mathcal T$ is called a tilting subcategory if the following conditions are satisfied:
\begin{itemize}
\item[$(1)$] $\emph{pd}_{\mathcal{E}}(\mathcal T)\leq n;$

\item[$(2)$] $\mathcal T$ is an generator for $\mathcal T_{\mathcal{E}}^{\bot}.$

Where the definition of $\emph{pd}_{\mathcal{E}}(\mathcal T), \ \mathcal T_{\mathcal{E}}^{\bot}, \ \mathcal{P}(\mathcal{E})$ and $\check{\mathcal T_{n}}$ are the same as that in \cite{be} and Section 3 in this paper.
\end{itemize}
\end{definition}

Applying Theorem \ref{baz} and Theorem \ref{arc}, we have the following Bazzoni type characterization of tilting subcategories and Auslander-Reiten correspondence in this situation.

\begin{theorem}\label{ptbaz}
Let $(\C,Ph_{\mathcal{E}},\s)$  be as above with enough projectives and injectives. Assume that  $\mathcal T\subset\C$ such that $\mathcal T$ is closed under direct summands and every object in $\mathcal T_{\mathcal{E}}^{\bot}$ has a right $\mathcal T-$approximation. Then $\mathcal T$ is a tilting subcategory iff $Pres^{n}(\mathcal T)=\mathcal T_{\mathcal{E}}^{\bot}.$
\end{theorem}

\begin{theorem}\label{ptarc}
Let $(\C,Ph_{\mathcal{E}},\s)$  be as above with enough projectives and injectives. The assignments $\mathcal T\mapsto \mathcal T^{\bot}$ and $\X\mapsto \sideset{^{\bot}}{}{\mathop{\X}}\cap\X$ induce a one-to-one correspondence between the class of tilting subcategories and coresolving covariantly finite subcategories $\X,$ which are closed under summands with $\check{\X_{n}}=\C.$
\end{theorem}

Note that one can also define cotilting subcategories in this setting, and the dual results for cotilting subcategories also hold for $(\C,Ph_{\mathcal{E}},\s).$

Bin Zhu\\
Department of Mathematical Sciences, Tsinghua University,
100084 Beijing, P. R. China.\\
E-mail: \verb"zhu-b@mail.tsinghua.edu.cn"\\[0.3cm]
Xiao Zhuang\\
Department of Mathematical Sciences, Tsinghua University,
100084 Beijing, P. R. China.\\
E-mail: \verb"zhuangx16@mails.tsinghua.edu.cn"

\end{document}